\documentclass[11pt,reqno]{amsart}
\usepackage{amssymb,amscd,amsbsy,mathrsfs}
\newcommand{\lb}{\linebreak}

\renewcommand{\a}{\alpha}
\renewcommand{\b}{\beta}
\newcommand{\g}{\gamma}

\newcommand{\e}{\varepsilon}
\newcommand{\vk}{\varkappa}
\newcommand{\z}{\zeta}

\renewcommand{\l}{\lambda}

\newcommand{\s}{\sigma}
\renewcommand{\t}{\tau}

\newcommand{\f}{\varphi}
\renewcommand{\o}{\omega}

\newcommand{\D}{\Delta}
\renewcommand{\L}{\Lambda}

\renewcommand{\O}{\Omega}

\newcommand{\B}{{\mathscr B}}

\newcommand{\E}{{\mathscr E}}

\newcommand{\F}{{\mathscr F}}

\newcommand{\h}{{\mathscr H}}

\newcommand{\X}{{\mathscr X}}
\newcommand{\Y}{{\mathscr Y}}

\newcommand{\W}{{\mathscr W}}

\newcommand{\C}{{\Bbb C}}
\newcommand{\T}{{\Bbb T}}

\newcommand{\R}{{\Bbb R}}
\newcommand{\Z}{{\Bbb Z}}

\newcommand{\0}{{\boldsymbol{0}}}

\newcommand{\bs}{\boldsymbol}

\newcommand{\bS}{{\boldsymbol S}}

\newcommand{\rf}[1]{(\ref{#1})}

\newcommand{\df}{\stackrel{\mathrm{def}}{=}}

\newcommand{\supp}{\operatorname{supp}}

\newcommand{\trace}{\operatorname{trace}}
\newcommand{\rank}{\operatorname{rank}}
\newcommand{\const}{\operatorname{const}}

\newcommand{\eeq}{\end{equation}}
\newcommand{\beq}{\begin{equation}}
\newcommand{\bay}{\begin{eqnarray}}
\newcommand{\ba}{\begin{align*}}
\newcommand{\ea}{\end{align*}}
\newcommand{\ey}{\end{eqnarray}}
\newcommand{\bey}{\begin{eqnarray*}}
\newcommand{\eey}{\end{eqnarray*}}

\newcommand{\be}{\infty}

\newcommand{\bl}{\blacksquare}
\newcommand{\ess}{\operatorname{ess}}

\newcommand{\Pf}{{\bf Proof. }}
\newcommand{\im}{\operatorname{Im}}

\newcommand{\ov}{\overline}

\newtheorem{thm}{\hspace{\parindent}Theorem}[section]

\newtheorem{cor}[thm]{\hspace{\parindent}Corollary}
\newtheorem{lem}[thm]{\hspace{\parindent}Lemma}

\usepackage{amsmath}

\makeatletter
\def\upintkern@{\mkern-7mu\mathchoice{\mkern-3.5mu}{}{}{}}
\def\upintdots@{\mathchoice{\mkern-4mu\@cdots\mkern-4mu}%
 {{\cdotp}\mkern1.5mu{\cdotp}\mkern1.5mu{\cdotp}}%
 {{\cdotp}\mkern1mu{\cdotp}\mkern1mu{\cdotp}}%
 {{\cdotp}\mkern1mu{\cdotp}\mkern1mu{\cdotp}}}

\newcommand{\UpMultiIntegral}[1]{%
  \edef\ints@c{\noexpand\upintop
    \ifnum#1=\z@\noexpand\upintdots@\else\noexpand\upintkern@\fi
    \ifnum#1>\tw@\noexpand\upintop\noexpand\upintkern@\fi
    \ifnum#1>\thr@@\noexpand\upintop\noexpand\upintkern@\fi
    \noexpand\upintop
    \noexpand\ilimits@
  }%
  \futurelet\@let@token\ints@a
}
\makeatother

\DeclareFontFamily{OMX}{mdbch}{}
\DeclareFontShape{OMX}{mdbch}{m}{n}{ <->s * [0.8]  mdbchr7v }{}
\DeclareFontShape{OMX}{mdbch}{b}{n}{ <->s * [0.8]  mdbchb7v }{}
\DeclareFontShape{OMX}{mdbch}{bx}{n}{<->ssub * mdbch/b/n}{}

\DeclareSymbolFont{uplargesymbols}{OMX}{mdbch}{m}{n}
\SetSymbolFont{uplargesymbols}{bold}{OMX}{mdbch}{b}{n}
\DeclareMathSymbol{\upintop}{\mathop}{uplargesymbols}{82}
\DeclareMathSymbol{\upointop}{\mathop}{uplargesymbols}{"48}

\DeclareFontEncoding{MDB}{}{}
\DeclareFontFamily{MDB}{mdbch}{}
\DeclareFontShape{MDB}{mdbch}{m}{n}{ <->s * [0.8]  mdbchrmb }{}
\DeclareFontShape{MDB}{mdbch}{b}{n}{ <->s * [0.8]  mdbchbmb }{}
\DeclareFontShape{MDB}{mdbch}{bx}{n}{<->ssub * mdbch/b/n}{}
\DeclareFontSubstitution{MDB}{cmr}{m}{n}
\DeclareSymbolFont{mathdesignB}{MDB}{mdbch}{m}{n}%
\SetSymbolFont{mathdesignB}{bold}{MDB}{mdbch}{b}{n}%
\DeclareMathSymbol{\upintclockwise}{\mathop}{mathdesignB}{128}
\DeclareMathSymbol{\upointclockwise}{\mathop}{mathdesignB}{130}
\DeclareMathSymbol{\upointctrclockwise}{\mathop}{mathdesignB}{132}
\DeclareMathSymbol{\upoiint}{\mathop}{mathdesignB}{134}
\DeclareMathSymbol{\upoiiint}{\mathop}{mathdesignB}{136}

\makeatletter
\newcommand{\upint}{\DOTSI\upintop\ilimits@}
\newcommand{\upoint}{\DOTSI\upointop\ilimits@}
\makeatother

\theoremstyle{remark}

\newtheorem*{rem*}{Remark}

\newcommand\fM{\frak M}

\newcommand\dg{\frak D}

\newcommand{\fI}{{\frak I}}

\newcommand\sG{\mathscr G}

\newcommand\mB{\mathcal{B}}

\newcommand{\Lip}{\operatorname{Lip}}

\begin{document}

\newcommand{\vse}{\vspace{.2in}}
\numberwithin{equation}{section}

\title{Functions of noncommuting self-adjoint operators under perturbation
and estimates of triple operator integrals}
\author{A.B. Aleksandrov, F.L. Nazarov and V.V. Peller}
\thanks{The first author is partially supported by RFBR grant 14-01-00198;
the second author is partially supported
by NSF grant DMS 1265623; the third author is partially supported by NSF grant DMS 1300924}
\thanks{Corresponding author: V.V. Peller; email: peller@math.msu.edu}

\begin{abstract}
We define functions of noncommuting self-adjoint operators with the help of double operator integrals. We are studying the problem to find conditions on a function $f$ on $\R^2$, for which the map $(A,B)\mapsto f(A,B)$ is Lipschitz in the operator norm and in Schatten--von Neumann norms $\bS_p$. It turns out that for functions $f$ in the Besov class $B_{\be,1}^1(\R^2)$, the above map is Lipschitz in the $\bS_p$ norm for $p\in[1,2]$. However, it is not Lipschitz in the operator norm, nor in the $\bS_p$ norm for $p>2$. The main tool is triple operator integrals. To obtain the results, we introduce new Haagerup-like tensor products of $L^\be$ spaces and 
obtain Schatten--von Neumann norm estimates of triple operator integrals. We also obtain similar results for functions of noncommuting unitary operators.
\end{abstract}

\maketitle

\

\begin{center}
{\Large Contents}
\end{center}

\

\begin{enumerate}
\item[1.] Introduction \quad\dotfill \pageref{In}
\item[2.] Preliminaries \quad\dotfill \pageref{prel}
\item[3.] Triple operator integrals
\quad\dotfill \pageref{Stoi}
\item[4.] Schatten--von Neumann properties of triple operator integrals  \quad\dotfill \pageref{SvNtoi}
\item[5.] Haagerup like tensor products and triple operator integrals \quad\dotfill \pageref{Ttoi}
\item[6.] When do the divided differences $\dg^{[1]}f$ and 
$\dg^{[2]}f$ belong to Haagerup-like tensor products? \quad\dotfill \pageref{ddiff}
\item[7.] Lipschitz type estimates in the case $1\le p\le2$ \quad\dotfill \pageref{ple2}
\item[8.] No Lipschitz type estimates in the operator norm and in the 
$\bS_p$ norm for \lb$p>2\,\,$! \quad\dotfill \pageref{Bp>2}
\item[9.] Two counterexamples \quad\dotfill \pageref{2c}
\item[10.] Points of local Lipschitzness \quad\dotfill \pageref{tlL}
\item[11.] A sufficient condition for Lipschitz type estimates \quad\dotfill \pageref{dostu}
\item[12.] Functions of noncommuting unitary operators \quad\dotfill \pageref{unit}
\item[] References \quad\dotfill \pageref{bibl}
\end{enumerate}

\

\setcounter{section}{0}
\section{\bf Introduction}
\setcounter{equation}{0}
\label{In}

\medskip

In this paper we study the behavior of functions $f(A,B)$ of (not necessarily commuting) self-adjoint operators $A$ and $B$ under perturbations of the pair 
$(A,B)$. 

In the case of commuting self-adjoint operators $A$ and $B$, one can define functions $f(A,B)$ for all bounded Borel on $\R^2$ (actually, it suffices to have a function $f$ defined on the joint spectrum $\s(A,B)$ of $A$ and $B$, which is a subset of the cartesian product $\s(A)\times\s(B)$ of the spectra of 
$A$ and $B$) by the formula
$$
f(A,B)\df\iint f(x,y)\,dE_{A,B}(x,y),
$$
where $E_{A,B}$ is the joint spectral measure of $A$ and $B$ (see \cite{BS0}).

If the self-adjoint operators $A$ and $B$ do not commute, we define the function
$f(A,B)$ as the {\it double operator integral}
\bay
\label{oprdv}
f(A,B)\df\iint f(x,y)\,dE_A(x)\,dE_B(y),
\ey
where $E_A$ and $E_B$ are the spectral measures of $A$ and $B$. The theory of double operator integrals was developed by Birman and Solomyak in \cite{BS1}, \cite{BS2}, and \cite{BS3}. Unlike in the case of commuting operators, the functions $f(A,B)$ cannot be defined for arbitrary bounded Borel functions $f$. For the integral in \rf{oprdv} to make sense, the function $f$ has to be a {\it Schur multiplier}. In \S~\ref{prel} we give a brief introduction in double operator integrals and Schur multipliers.

Let us briefly summarize known results on perturbations of functions of one operator and functions of commuting operators. A function $f$ on the real line $\R$ is called {\it operator Lipschitz} if 
$$
\|f(A)-f(B)\|\le\const\|A-B\|
$$
for arbitrary self-adjoint operators $A$ and $B$ (bounded or, possibly, unbounded).
It was shown in \cite{F} that a Lipschitz function $f$ (i.e., a function satisfying
the inequality $|f(x)-f(y)|\le\const|x-y|$, $x,\,y\in\R$) does not have to be operator Lipschitz. It turned out later (see \cite{Mc} and \cite{K}) that the function $x\mapsto|x|$ is not operator Lipschitz. Note also that in \cite{Pe2} and \cite{Pe3} necessary conditions for $f$ to be operator Lipschitz are found that are based on the trace class description of Hankel operators (see \cite{Pe1} and \cite{Pe5}).

Among various sufficient conditions for operator Lipschitzness we mention
the one found in \cite{Pe2} (see also \cite{Pe3}) in terms of Besov spaces:
if $f$ belongs to the Besov class $B_{\be,1}^1(\R)$, then $f$ is operator Lipschitz (see \S~\ref{prel} for a brief introduction in Besov spaces).

It is well known that $f$ is operator Lipschitz if and only if it possesses the property
$$
A-B\in\bS_1\quad\Longrightarrow\quad f(A)-f(B)\in\bS_1.
$$
Moreover, if $f$ is operator Lipschitz, then it is also trace norm Lipschitz, i.e.,
$$
\|f(A)-f(B)\|_{\bS_1}\le\const\|A-B\|_{\bS_1}.
$$
Here $\bS_1$ is trace class. We are going to use the notation $\bS_p$ for Schatten--von Neumann classes and we refer the reader to \cite{GK} for detailed information about such ideals of operators.

If we consider Lipschitz type estimates in the Schatten--von Neumann norm 
$\bS_p$, \lb$1<p<\be$, the situation is different. A classical result (see \cite{BS3})
says that if $A-B$ belongs to the Hilbert-Schmidt class $\bS_2$ and $f$ is a Lipschitz function, than
$$
\|f(A)-f(B)\|_{\bS_2}\le\|f\|_{\Lip}\|A-B\|_{\bS_2}.
$$
Recently it was shown in \cite{PoS} that such a Lipschitz type estimate also holds in the $\bS_p$ norm for $p\in(1,\be)$ with a constant on the right-hand side that depends on $p$.

It turns out, however, that the situation is entirely different if we proceed from Lipschitz functions to H\"older functions. It was shown in \cite{AP2}  that if $f$ belongs to the H\"older class $\L_\a(\R)$, $0<\a<1$, i.e., $|f(x)-f(y)|\le\const|x-y|^\a$, $x,\,y\in\R$, then $f$ is necessarily {\it operator H\"older of order} $\a$, i.e.,
$$
\|f(A)-f(B)\|\le\const\|A-B\|^\a
$$
for all self-adjoint operators $A$ and $B$ on Hilbert space with bounded $A-B$. 
Note that in \cite{AP2} sharp results were also obtained for functions in the space $\L_\o$ for an arbitrary modulus of continuity $\o$.
Similar (slightly weaker) results were obtained independently in \cite{FN}.  

It was proved in \cite{AP3}  that for $f\in\L_\a(\R)$, $0<\a<1$, $p>1$, and for self-adjoint operators $A$ and $B$ with $A-B\in\bS_p$, the operator $f(A)-f(B)$ must be in $\bS_{p/\a}$ and the following inequality holds:
$$
\|f(A)-f(B)\|_{\bS_{p/\a}}\le\const\|f\|_{\L_\a}\|A-B\|_{\bS_p}^\a.
$$
Let us also mention that in \cite{AP3} more general results for operator ideals were obtained as well.

It turns out that the situation for functions of normal operators or, which is the equivalent, for functions of two commuting self-adjoint operators is more complicated and requires different techniques. Nevertheless, in \cite{APPS}  analogs of the above mentioned results were obtained for normal operators and functions on the plane. In particular, it was shown in \cite{APPS} that if $f$ belongs to the Besov class $B_{\be,1}^1(\R^2)$, then
$$
\|f(N_1)-f(N_2)\|\le\const\|f\|_{B_{\be,1}^1}\|N_1-N_2\|
$$
for arbitrary normal operators $N_1$ and $N_2$.

However, the methods of \cite{APPS} do not work in the case of functions of more than two commuting self-adjoint operators. New methods were found in \cite{NP} to obtain analogs of the above results for functions of $n$-tuples of commuting self-adjoint operators.

Note also that the Lipschitz type estimate for Lipschitz functions in the norm of 
$\bS_p$ with $1<p<\be$ was generalized in \cite{KPSS} to $n$-tuples of commuting self-adjoint operators.

The subject of this paper is estimates of $f(A_1,B_1)-f(A_2,B_2)$, where 
$(A_1,B_1)$, and $(A_2,B_2)$ are pairs of (not necessarily commuting) self-adjoint operators. Here we consider the pair $(A_2,B_2)$ as a perturbation of the pair $(A_1,B_1)$.

The main tool is estimates of triple operator integrals 
(see \S~\ref{Stoi} and \S~\ref{Ttoi} for a detailed discussion of triple operator integrals). To establish a Lipschitz type inequality in trace norm
$$
\|f(A_1,B_1)-f(A_2,B_2)\|_{\bS_1}
\le\const\max\{\|A_1-A_2\|_{\bS_1},\|A_1-A_2\|_{\bS_1}\}
$$
for functions in the Besov class $B_{\be,1}^1(\R^2)$, we would like to apply the following formula
\begin{align}
\label{predst}
f(A_1,B_1)&-f(A_2,B_2)=
\iiint\big(\dg^{[1]}f\big)(x_1,x_2,y)\,dE_{A_1}(x_1)(A_1-A_2)\,dE_{A_2}(x_2)\,dE_{B_1}(y)
\nonumber\\[.2cm]
&+\iiint\big(\dg^{[2]}f\big)(x,y_1,y_2)\,dE_{A_2}(x)\,dE_{B_1}(y_1)(B_1-B_2)\,dE_{B_2}(y_2).
\end{align}
Here the divided differences $\dg^{[1]}f$ and $\dg^{[2]}f$ are defined by
$$
\dg^{[1]}f(x_1,x_2,y)\df\frac{f(x_1,y)-f(x_2,y)}{x_1-x_2}
\quad\mbox{and}\quad
\dg^{[2]}f(x,y_1,y_2)\df\frac{f(x,y_1)-f(x,y_2)}{y_1-y_2},
$$
Triple operator integrals can be defined when we integrate functions in the projective tensor product $L^\be\hat\otimes L^\be\hat\otimes L^\be$ or in the Haagerup tensor product $L^\be\!\otimes_{\rm h}\!L^\be\!\otimes_{\rm h}\!L^\be$
(see the definitions in \S~\ref{Stoi}).
It turns out the for $f\in B_{\be,1}^1(\R^2)$, the divided differences 
$\dg^{[1]}f$ and $\dg^{[2]}f$ do not have to belong to the Haagerup tensor product $L^\be\!\otimes_{\rm h}\!L^\be\!\otimes_{\rm h}\!L^\be$ (and a fortiori to the
projective tensor product $L^\be\hat\otimes L^\be\hat\otimes L^\be$). This will be proved in \S~\ref{2c}.

To overcome the problem, we introduce in \S~\ref{Ttoi} Haagerup-like tensor products of the first kind and of the second kind, define triple operator products for functions in such Haagerup like tensor products, and prove in \S~\ref{ddiff} that  that for $f\in B_{\be,1}^1(\R^2)$, the divided difference
$\dg^{[1]}f$ belongs 
to the Haagerup-like tensor product of the first kind, while $\dg^{[2]}f$ belongs 
to the Haagerup-like tensor product of the second kind.

We obtain in \S~\ref{ple2} the following Lipschitz type inequality
\bay
\|f(A_1,B_1)-f(A_2,B_2)\|_{\bS_p}
\le\const\|f\|_{B_{\be,1}^1}\max\{\|A_1-A_2\|_{\bS_p},\|A_1-A_2\|_{\bS_p}\}
\ey
for $p\in[1,2]$. To prove this inequality, we obtain in \S~\ref{SvNtoi} certain Schatten--von Neumann estimates for triple operator integrals
\bay
\label{tPsi}
\iiint\Psi(x_1,x_2,x_3)\,dE_1(x_1)T\,dE_2(x_2)R\,dE_3(x_3).
\ey
In particular, we show in \S~\ref{SvNtoi} that if $\Psi$ belongs to the Haagerup tensor product \lb$L^\be\!\otimes_{\rm h}\!L^\be\!\otimes_{\rm h}\!L^\be$,
$T$ is a bounded operator and $R\in\bS_p$ with $p\ge2$, then the triple operator integral \rf{tPsi} belongs to $\bS_p$. However, for $p<2$ this is not true which will be proved in \S~\ref{2c}. We also establish in \S~\ref{SvNtoi} that if
$\Psi\in L^\be\!\otimes_{\rm h}\!L^\be\!\otimes_{\rm h}\!L^\be$, $T\in\bS_p$, $R\in\bS_q$, and $1/p+1/q\le1/2$, then the triple operator integral \rf{tPsi} belongs to $\bS_r$, where $1/r=1/p+1/q$.

In \S~\ref{Bp>2} we show that a Lipschitz type inequality in the norm of $\bS_p$ with $p>2$ does not hold. The same is true in the operator norm.

It turns out, however, that in the operator norm (as well as in the $S_p$ norm for any $p>0$) there are points of Lipschitzness of the map $(A,B)\mapsto f(A,B)$ for $f\in B_{\be,1}^1(\R^2)$. We prove in \S~\ref{tlL} that the pairs $(\a I, \b I)$ with
$\a,\,\b\in\R$ are points of Lipschitzness.

We find in \S~\ref{dostu} a sufficient condition on a function $f$ under which the Lipschitz type inequality in the operator norm (as well as in the norms of $\bS_p$ with $p\ge1$) holds.

Finally, in \S~\ref{unit} we obtain similar results for functions of noncommuting unitary operators.

In \S~\ref{prel} we collect necessary information on Besov classes, integration of vector functions with respect to spectral measures, double operator integrals, and 
functions of noncommuting operators.

\

\section{\bf Preliminaries}
\setcounter{equation}{0}
\label{prel}

\

In this section we collect necessary information on function spaces, operator ideals, and double operator integrals. 

\medskip

{\bf 2.1.1. Besov classes of functions on Euclidean spaces and Littlewood--Paley type expansions.} The technique of Littlewood--Paley type expansions of functions or distributions on Euclidean spaces 
is a very important tool in Harmonic Analysis. 

Let $w$ be an infinitely differentiable function on $\R$ such
that
\bay
\label{w}
w\ge0,\quad\supp w\subset\left[\frac12,2\right],\quad\mbox{and} \quad w(s)=1-w\left(\frac s2\right)\quad\mbox{for}\quad s\in[1,2].
\ey

We define the functions $W_n$, $n\in\Z$, on $\R^d$ by 
$$
\big(\F W_n\big)(x)=w\left(\frac{\|x\|_2}{2^n}\right),\quad n\in\Z, \quad x=(x_1,\cdots,x_d),
\quad\|x\|_2\df\left(\sum_{j=1}^dx_j^2\right)^{1/2},
$$
where $\F$ is the {\it Fourier transform} defined on $L^1\big(\R^n\big)$ by
$$
\big(\F f\big)(t)=\!\int\limits_{\R^n} f(x)e^{-{\rm i}(x,t)}\,dx,\!\quad 
x=(x_1,\cdots,x_d),
\quad t=(t_1,\cdots,t_d), \!\quad(x,t)\df \sum_{j=1}^dx_jt_j.
$$
Clearly,
$$
\sum_{n\in\Z}(\F W_n)(t)=1,\quad t\in\R^d\setminus\{0\}.
$$

With each tempered distribution $f\in{\mathscr S}^\prime\big(\R^d\big)$, we
associate the sequence $\{f_n\}_{n\in\Z}$,
\bay
\label{fn}
f_n\df f*W_n.
\ey
The formal series
$
\sum_{n\in\Z}f_n
$
is a Littlewood--Paley type expansion of $f$. This series does not necessarily converge to $f$. Note that in this paper we mostly deal with Besov spaces $B_{\be,1}^1(\R^d)$.
For functions $f$ in $B_{\be,1}^1(\R^d)$,
$$
f(x)-f(y)=\sum_{n\in\Z}\big(f_n(x)-f_n(y)\big),\quad x,~y\in\R^d,
$$
and the series on the right converges uniformly.

Initially we define the (homogeneous) Besov class $\dot B^s_{p,q}\big(\R^d\big)$,
$s>0$, $1\le p,\,q\le\be$, as the space of all
$f\in{\mathscr S}^\prime(\R^n)$
such that
\bay
\label{Wn}
\{2^{ns}\|f_n\|_{L^p}\}_{n\in\Z}\in\ell^q(\Z)
\ey
and put
$$
\|f\|_{B^s_{p,q}}\df\big\|\{2^{ns}\|f_n\|_{L^p}\}_{n\in\Z}\big\|_{\ell^q(\Z)}.
$$
According to this definition, the space $\dot B^s_{p,q}(\R^n)$ contains all polynomials
and all polynomials $f$ satisfy the equality $\|f\|_{B^s_{p,q}}=0$. Moreover, the distribution $f$ is determined by the sequence $\{f_n\}_{n\in\Z}$
uniquely up to a polynomial. It is easy to see that the series 
$\sum_{n\ge0}f_n$ converges in ${\mathscr S}^\prime(\R^d)$.
However, the series $\sum_{n<0}f_n$ can diverge in general. It can easily be proved that the series
\bay
\label{ryad}
\sum_{n<0}\frac{\partial^r f_n}{\partial x_1^{r_1}\cdots\partial x_d^{r_d}},\qquad \mbox{where}\quad r_j\ge0,\quad\mbox{for}\quad
1\le j\le d,\quad\sum_{j=1}^dr_j=r,
\ey
converges uniformly on $\R^d$ for every nonnegative integer
$r>s-d/p$. Note that in the case $q=1$ the series \rf{ryad}
converges uniformly, whenever $r\ge s-d/p$.

Now we can define the modified (homogeneous) Besov class $B^s_{p,q}\big(\R^d\big)$. We say that a distribution $f$
belongs to $B^s_{p,q}(\R^d)$ if \rf{Wn} holds and
$$
\frac{\partial^r f}{\partial x_1^{r_1}\cdots\partial x_d^{r_d}}
=\sum_{n\in\Z}\frac{\partial^r f_n}{\partial x_1^{r_1}\cdots\partial x_d^{r_d}},\quad
\mbox{whenever}\quad 
r_j\ge0,\quad\mbox{for}\quad
1\le j\le d,\quad\sum_{j=1}^dr_j=r.
$$
in the space ${\mathscr S}^\prime\big(\R^d\big)$, where $r$ is
the minimal nonnegative integer such that $r>s-d/p$ ($r\ge s-d/p$ if $q=1$). Now the function $f$ is determined uniquely by the sequence $\{f_n\}_{n\in\Z}$ up
to a polynomial of degree less than $r$, and a polynomial $g$ belongs to 
$B^s_{p,q}\big(\R^d\big)$
if and only if $\deg g<r$.

As we have already mentioned, in this paper we deal with Besov classes 
$B_{\be,1}^1(\R^d)$. They can also be defined in the following way:

Let $X$ be the set of all continuous functions $f\in L^\be(\R^d)$ such that $|f|\le1$ and 
$\supp\F f\subset\{\xi\in\R^d:~\|\xi\|\le1\}$.
Then 
$$
B^1_{\be1}(\R^d)=\left\{c+\sum_{n=1}^\be \a_n\s_n^{-1}(f_n(\s_nx)-f(0)):
~c\in\C,~f_n\in X, ~\s_n>0,~\sum_{n=1}^\be|\a_n|<\be\right\}.
$$

Note that the functions $f_\s$, $f_\s(x)=f(\s x)$, $x\in\R^d$, have the following properties: \lb$f_\s\in L^\be(\R^d)$ and 
$\supp\F f\subset\{\xi\in\R^d:~\|\xi\|\le\s\}$. Such functions can be characterized by the following Paley--Wiener--Schwartz type theorem  (see \cite{R}, Theorem 7.23 and exercise 15 of Chapter 7):

{\it Let $f$ be a continuous function
on $\R^d$ and let $M,\,\s>0$. The following statements are equivalent:

{\em(i)} $|f|\le M$ and $\supp\F f\subset\{\xi\in\R^d:\|\xi\|_2\le\s\}$;

{\em(ii)} $f$ is a restriction to $\R^d$ of an entire function on $\C^d$ such that 
$$
|f(z)|\le Me^{\s\|\im z\|_2}
$$
for all $z\in\C^d$.}


Besov classes admit many other descriptions.
We give here the definition in terms of finite differences.
For $h\in\R^d$, we define the difference operator $\D_h$,
$$
(\D_hf)(x)=f(x+h)-f(x),\quad x\in\R^d.
$$
It is easy to see that $B_{p,q}^s\big(\R^d\big)\subset L^1_{\rm loc}\big(\R^d\big)$ for every $s>0$
and $B_{p,q}^s\big(\R^d\big)\subset C\big(\R^d\big)$ for every $s>d/p$. Let $s>0$ and let $m$ be the integer such that $m-1\le s<m$.
The Besov space $B_{p,q}^s\big(\R^d\big)$ can be defined as the set of
functions $f\in L^1_{\rm loc}\big(\R^d\big)$ such that
$$
\int_{\R^d}|h|^{-d-sq}\|\D^m_h f\|_{L^p}^q\,dh<\be\quad\mbox{for}\quad q<\be
$$
and
$$
\sup_{h\not=0}\frac{\|\D^m_h f\|_{L^p}}{|h|^s}<\be\quad\mbox{for}\quad q=\be.
$$
However, with this definition the Besov space can contain polynomials of higher degree than in the case of the first definition given above.

We refer the reader to \cite{Pee} and \cite{Tr} for more detailed information on Besov spaces.

\medskip

{\bf 2.1.2. Besov classes of periodic functions.} Studying periodic functions on $\R^d$ is equivalent to studying functions on the $d$-dimensional torus $\T^d$. To define Besov spaces on $\T^d$, we consider a function $w$ satisfying \rf{w} and define the trigonometric polynomials $W_n$, $n\ge0$, by
$$
W_n(\z)\df\sum_{j\in\Z^d}w\left(\frac{\|\z\|_2}{2^n}\right)\z^j,\quad n\ge1,
\quad W_0(\z)\df\sum_{\{j:\|j\|_2\le1\}}\z^j,
$$
where 
$$
\z=(\z_1,\cdots,\z_d)\in\T^d,\quad j=(j_1,\cdots,j_d),\quad\mbox{and}\quad
\|\z\|_2=\big(|\z_1|^2+\cdots+|\z_d|^2\big)^{1/2}.
$$
For a distribution $f$ on $\T^d$ we put
$$
f_n=f*W_n,\quad n\ge0,
$$
and we say that $f$ belongs the Besov class $B_{p,q}^s(\T^d)$, $s>0$, 
$1\le p,\,q\le\be$, if
\bay
\label{Bperf}
\big\{2^ns\|f_n\|_{L^p}\big\}_{n\ge0}\in\ell^q.
\ey

Note that locally the Besov space $B_{p,q}^s(\R^d)$ coincides with the Besov space
$B_{p,q}^s$ of periodic functions on $\R^d$.

\medskip

{\bf 2.2. Integration of vector functions with respect to spectral measures.} Let $E$ be a spectral measure on a $\s$-algebra of subsets of $\O$ that takes values in the set of orthogonal projections on a Hilbert space $\h$. It is well known that for a scalar function $f$ in $L^\be(E)$ the integral $\int f\,dE$ admits the estimate
$$
\left\|\int_\O f(\o)\,dE(\o)\right\|\le\|f\|_{L^\be(E)}.
$$
We would like to be able to integrate $\h$-valued functions to get vectors in $\h$. However, it is easy to see that unlike the case of scalar functions it is impossible to define an integral of an arbitrary bounded measurable $\h$-valued function.
We consider the {\it projective tensor product} $L^\be(E)\hat\otimes\h$, which consists of 
$\h$-valued functions $f$ that admit a representation of the form
\bay
\label{EtH}
f(\o)=\sum_n \f_n(\o)v_n,\quad\o\in\O,
\ey
where $\f_n\in L^\be(E)$, $v_n\in\h$, and
\bay
\label{pro}
\sum_n\|\f_n\|_{L^\be(E)}\|v_n\|_\h<\be.
\ey
The norm of $f$ in $L^\be(E)\hat\otimes\h$ is defined as the infimum of the left-hand side of \rf{pro} over all representations of the form \rf{EtH}. For an $\h$-valued function $f$ of the form \rf{EtH}, we put
\bay
\label{SiE}
\int_\O\big(dE(\o)f(\o)\big)\df\sum_n\left(\int_\O\f_n(\o)\,dE(\o)\right)v_n.
\ey
It follows from \rf{pro} that the series on the right-hand side of \rf{SiE} converges absolutely in the norm of $\h$. Let us show that the integral is well defined.

\begin{thm}
\label{iHvf}
The right-hand side of {\em\rf{SiE}} does not depend on the choice of a representation of $f$ of the form {\em\rf{EtH}}.
\end{thm}

\Pf Clearly, it suffices to prove that if
\bay
\label{tozh}
\sum_n \f_n(\o)v_n=\0,\quad\o\in\O,
\ey
and \rf{pro} holds, then
$$
\sum_n\left(\int_\O\f_n(\o)\,dE(\o)\right)v_n=\0.
$$

Without loss of generality we can assume that the functions $\f_n$ are defined everywhere,
$\|\f_n\|_{L^\be(E)}=\sup|\f_n|$ for all $n$ and equality \rf{tozh} holds for all $\o$ in $\O$.

Consider first the special case when
the range of the vector function $\f=\{\f_n\}$ is
finite. Let $\big\{\l(k)=\{\l_n(k)\}: k=1,2,\dots,N\big\}$ be the set of values of $\f$. Then  
$\sum_{n=1}^\be\l_n(k)v_n=0$ for all $k$.
Put $P_k\df E\{\o:~\f(\o)=\l_k\}$. We have to prove that
$$
\sum_{n=1}^\be\left(\sum_{k=1}^N \l_n(k)P_k\right)v_n=0.
$$
We have
$$
\sum_{n=1}^\be\left(\sum_{k=1}^N \l_n(k)P_k\right)v_n=
\sum_{k=1}^N\left(\sum_{n=1}^\be \l_n(k)P_k v_n\right)=
\sum_{k=1}^N P_k\left(\sum_{n=1}^\be \l_n(k)v_n\right)=0.
$$

We proceed now to the general case. 

Clearly, we can construct a sequence $\{\f^{(j)}\}$
of vector functions such that the range of $\f^{(j)}$ is a finite subset of
the range of $\f$ for all $j$ and $\big|\f_n-\f^{(j)}_n\big|\le2^{-j}\|\f_n\|_{L^\be(E)}$ everywhere for $n=1,2,\dots,j$.
We have
\begin{align*}
\sum\limits_{n=1}^\be\left(\int\f_n\,dE\right)v_n&=
\sum\limits_{n=1}^\be\left(\int\big(\f_n-\f_n^{(j)}\big)\,dE\right)v_n\\[.2cm]
&+
\sum\limits_{n=1}^\be\left(\int\f_n^{(j)}\,dE\right)v_n
=
\sum\limits_{n=1}^\be\left(\int\big(\f_n-\f_n^{(j)}\big)\,dE\right)v_n
\end{align*}
which follows from the special case considered above. Hence, 
\begin{align*}
\left\|\sum_{n=1}^\be\left(\int\f_n\,dE\right)v_n\right\|_\h&=
\left\|\sum\limits_{n=1}^\be\left(\int\big(\f_n-\f_n^{(j)}\big)\,dE\right)v_n\right\|_\h
\\[.2cm]
&\le\sum_{n=1}^j\left\|\int\big(\f_n-\f_n^{(j)}\big)\,dE\right\|\cdot\|v_n\|_\h\\[.2cm]
&+\sum_{n=j+1}^\be\left\|\int\big(\f_n-\f_n^{(j)}\big)\,dE\right\|\cdot\|v_n\|_\h\\[.2cm]
&\le\frac1{2^j}\sum_{n=1}^j\|\f_n\|_{L^\be(E)}\|v_n\|_\h+2\sum_{n=j+1}^\be\|\f_n\|_{L^\be(E)}\|v_n\|_\h\\[.2cm]
&\le\frac1{2^j}\sum_{n=1}^\be\|\f_n\|_{L^\be(E)}\|v_n\|_\h+2\sum_{n=j+1}^\be\|\f_n\|_{L^\be(E)}\|v_n\|_\h\to0\\[.2cm]
&\mbox{as}\quad j\to+\be. \quad\bl
\end{align*}

\medskip

{\bf 2.3. Double operator integrals.}
In this subsection we give a brief introduction to double  operator integrals. Double operator integrals appeared in the paper \cite{DK} by Daletskii and S.G. Krein. Later the beautiful theory of double operator integrals was developed by Birman and Solomyak in \cite{BS1}, \cite{BS2}, and \cite{BS3}, see also their survey \cite{BS5}.

Let $(\X,E_1)$ and $(\Y,E_2)$ be spaces with spectral measures $E_1$ and $E_2$
on a Hilbert space $\h$. The idea of Birman and Solomyak is to define first
double operator integrals
\bay
\label{doi}
\int\limits_\X\int\limits_\Y\Phi(x,y)\,d E_1(x)T\,dE_2(y),
\ey
for bounded measurable functions $\Phi$ and operators $T$
of Hilbert Schmidt class $\bS_2$. Consider the spectral measure $\E$ whose values are orthogonal projections on the Hilbert space $\bS_2$, which is defined by
$$
\E(\L\times\D)T=E_1(\L)TE_2(\D),\quad T\in\bS_2,
$$
$\L$ and $\D$ being measurable subsets of $\X$ and $\Y$. It was shown in \cite{BS} that $\E$ extends to a spectral measure on
$\X\times\Y$. If $\Phi$ is a bounded measurable function on $\X\times\Y$, we define the double operator integral \rf{doi} by
$$
\int\limits_\X\int\limits_\Y\Phi(x,y)\,d E_1(x)T\,dE_2(y)\df
\left(\,\,\int\limits_{\X\times\Y}\Phi\,d\E\right)T.
$$
Clearly,
$$
\left\|\int\limits_\X\int\limits_\Y\Phi(x,y)\,dE_1(x)T\,dE_2(y)\right\|_{\bS_2}
\le\|\Phi\|_{L^\be}\|T\|_{\bS_2}.
$$
If
$$
\int\limits_\X\int\limits_\Y\Phi(x,y)\,d E_1(x)T\,dE_2(y)\in\bS_1
$$
for every $T\in\bS_1$, we say that $\Phi$ is a {\it Schur multiplier of $\bS_1$ associated with
the spectral measures $E_1$ and $E_2$}.

In this case the transformer
\bay
\label{tra}
T\mapsto\int\limits_{\Y}\int\limits_{\X}\Phi(x,y)\,d E_2(y)\,T\,dE_1(x),\quad T\in \bS_2,
\ey
extends by duality to a bounded linear transformer on the space of bounded linear operators on $\h$
and we say that the function $\Psi$ on $\Y\times\X$ defined by
$$
\Psi(y,x)=\Phi(x,y)
$$
is {\it a Schur multiplier (with respect to $E_2$ and $E_1$) of the space of bounded linear operators}.
We denote the space of such Schur multipliers by $\fM(E_2,E_1)$.
The norm of $\Psi$ in $\fM(E_2,E_1)$ is, by definition, the norm of the
transformer \rf{tra} on the space of bounded linear operators.

In \cite{BS3} it was shown that if $A$ and $B$ are self-adjoint operators (not necessarily bounded) such that $A-B$ is bounded
 and if $f$ is a continuously differentiable
function on $\R$ such that the divided difference $\dg f$,
$$
\big(\dg f\big)(x,y)=\frac{f(x)-f(y)}{x-y},
$$
is a Schur multiplier
with respect to the spectral measures of $A$ and $B$, then
$$
f(A)-f(B)=\iint\big(\dg f\big)(x,y)\,dE_{A}(x)(A-B)\,dE_B(y)
$$
and
$$
\|f(A)-f(B)\|\le\const\|\dg f\|_{\fM(E_A,E_{B})}\|A-B\|,
$$
i.e., {\it $f$ is an operator Lipschitz function}.

It was established in \cite{Pe2} (see also \cite{Pe3}) that if $f$ belongs to the Besov class $B_{\be,1}^1(\R)$, then the divided difference $\dg f\in \fM(E_1,E_2)$
for arbitrary Borel spectral $E_1$ and $E_2$, and so 
\bay
\label{BesSch}
\|f(A)-f(B)\|\le\const\|f\|_{B_{\be,1}^1}\|A-B\|
\ey
for arbitrary self-adjoint operators $A$ and $B$.

There are different characterizations of the space $\fM(E_1,E_2)$ of Schur multipliers, see \cite{Pe2} and \cite{Pi}. In particular, $\Phi\in\fM(E_1,E_2)$ if and only if $\Phi$ belongs to the {\it Haagerup tensor product} 
$L^\be(E_1)\!\otimes_{\rm h}\!L^\be(E_2)$ of the spaces $L^\be(E_1)$ and
$L^\be(E_2)$, i.e., $\Phi$ admits a representation
\bay
\label{Ffps}
\Phi(x,y)=\sum_{j\ge0}\f_j(x)\psi_j(y),
\ey
where 
$$
\{\f_j\}_{j\ge0}\in L^\be(\ell^2)\quad\mbox{and}\quad
\{\psi_j\}_{j\ge0}\in L^\be(\ell^2).
$$
For such functions $\Phi$ it is easy to verify that
\bay
\label{skh}
\int\limits_\X\int\limits_\Y\Phi(x,y)\,dE_1(x)T\,dE_2(y)=
\sum_{j\ge0}\left(\,\int\limits_\X\f_j\,dE_1\right)T\left(\,\int\limits_\Y\psi_j\,dE_2\right)
\ey
and the series on the right converges in the weak operator topology.

In this paper we need the following easily verifiable sufficient condition:

\medskip

{\it If a function $\Phi$ on $\X\times\Y$ belongs to the {\it projective tensor
product}
$L^\be(E_1)\hat\otimes L^\be(E_2)$ of $L^\be(E_1)$ and $L^\be(E_2)$ (i.e., $\Phi$ admits a representation of the form {\em\rf{Ffps}}
whith \lb$\f_j\in L^\be(E_1)$, $\psi_j\in L^\be(E_2)$, and
$$
\sum_{j\ge0}\|\f_j\|_{L^\be}\|\psi_j\|_{L^\be}<\be),
$$
then $\Phi\in\fM(E_1,E_2)$ and}
\bay
\label{dous}
\|\Phi\|_{\fM(E_1,E_2)}\le\sum_{j\ge0}\|\f_j\|_{L^\be}\|\psi_j\|_{L^\be}.
\ey

For such functions $\Phi$, formula \rf{skh} holds
and the series on the right-hand side of \rf{skh} converges absolutely in the norm.

\medskip

{\bf2.4. Functions of noncommuting self-adjoint operators.} Let $A$ and $B$ be self-adjoint operators on Hilbert space and let $E_A$ and $E_B$ be their spectral measures. Suppose that $f$ is a function of two variables that is defined at least on $\s(A)\times\s(B)$. As we have already mentioned in the introduction, if
 $f$ is a Schur multiplier with respect to the pair $(E_A,E_B)$, we define the function $f(A,B)$ of $A$ and $B$ by
\bay
\label{fAB}
f(A,B)\df\iint f(x,y)\,dE_A(x)\,dE_B(y).
\ey
Note that this functional calculus $f\mapsto f(A,B)$ is linear, but not multiplicative.

If we consider functions of bounded operators, without loss of generality we may deal with periodic functions with a sufficiently large period. Clearly, we can rescale the problem and assume that our functions are $2\pi$-periodic in each variable.

If $f$ is a trigonometric polynomial of degree $N$, we can represent $f$ in the form
$$
f(x,y)=\sum_{j=-N}^Ne^{{\rm i}jx}\left(\sum_{k=-N}^N\hat f(j,k)e^{{\rm i}ky}\right).
$$
Thus 
$f$ belongs to the projective tensor product $L^\be\hat\otimes L^\be$ and
$$
\|f\|_{L^\be\hat\otimes L^\be}\le\sum_{j=-N}^N\sup_y
\left|\sum_{k=-N}^N\hat f(j,k)e^{{\rm i}ky}\right|
\le(1+2N)\|f\|_{L^\be}
$$
It follows easily from \rf{Bperf} that every periodic function $f$ of Besov class $B_{\be1}^1$ of periodic functions belongs to 
$L^\be\hat\otimes L^\be$, and so the operator $f(A,B)$ is well defined by \rf{fAB}.

\

\section{\bf Triple operator integrals}
\setcounter{equation}{0}
\label{Stoi}

\

Multiple operator integrals were considered by several mathematicians, see \cite{Pa}, \cite{St}. However, those definitions required very strong restrictions on the classes of functions that can be integrated. In \cite{Pe4} multiple operator integrals were defined for functions that belong to the (integral) projective tensor product of $L^\be$ spaces. Later in \cite{JTT} multiple operator integrals were defined for Haagerup tensor products of $L^\be$ spaces.

In this paper we deal with triple operator integrals. 
We consider here both approaches given in \cite{Pe4} and \cite{JTT}.

Let $E_1$, $E_2$, and $E_3$ be spectral measures on Hilbert space and let $T$ and $R$ be bounded linear operators on Hilbert space. Triple operator integrals are expressions of the following form:
\bay
\label{troi}
\int\limits_{\X_1}\int\limits_{\X_2}\int\limits_{\X_3} 
\Psi(x_1,x_2,x_3)\,dE_1(x_1)T\,dE_2(x_2)R\,dE_3(x_3).
\ey
Such integrals make sense under certain assumptions on $\Psi$, $T$, and $R$. The function $\Psi$ will be called the {\it integrand} of the triple operator integral.

Recall that the {\it projective tensor product} 
$L^\be(E_1)\hat\otimes L^\be(E_2)\hat\otimes L^\be(E_3)$ can be defined as the class of function $\Psi$ of the form
\bay
\label{pred}
\Psi(x_1,x_2,x_3)=\sum_n\f_n(x_1)\psi_n(x_2)\chi_n(x_3)
\ey
such that
\bay
\label{norma}
\sum_n\|\f_n\|_{L^\be(E_1)}\|\psi_n\|_{L^\be(E_2)}\|\chi_n\|_{L^\be(E_3)}<\be.
\ey
The norm $\|\Psi\|_{L^\be\hat\otimes L^\be\hat\otimes L^\be}$ of $\Psi$ is, by definition, the infimum of the left-hand side of \rf{norma} over all representations of the form \rf{pred}.

For $\Psi\in L^\be(E_1)\hat\otimes L^\be(E_2)\hat\otimes L^\be(E_3)$ of the form 
\rf{pred} the triple operator integral \rf{troi} was defined in \cite{Pe4} by
\begin{align}
\label{otoi}
\iiint \Psi(x_1,x_2,x_3)&\,dE_1(x_1)T\,dE_2(x_2)R\,dE_3(x_3)\nonumber\\
=&\sum_n\left(\int\f_n\,dE_1\right)T\left(\int\psi_n\,dE_2\right)R\left(\int\chi_n\,dE_3\right).
\end{align}
Clearly, \rf{norma} implies that the series on the right converges absolutely in the norm. The right-hand side of \rf{otoi} does not depend on the choice of a representation of the form \rf{pred}. Clearly,
$$
\left\|\iiint \Psi(x_1,x_2,x_3)\,dE_1(x_1)T\,dE_2(x_2)R\,dE_3(x_3)\right\|
\le\|\Psi\|_{L^\be\hat\otimes L^\be\hat\otimes L^\be}\|T\|\cdot\|R\|.
$$
Note that for $\Psi\in L^\be(E_1)\hat\otimes L^\be(E_2)\hat\otimes L^\be(E_3)$, triple operator integrals have the following properties:
\bay
\label{bep}
T\in\mB(\h),\quad R\in\bS_p,\quad1\le p<\be,\quad\Longrightarrow\quad 
\iiint\Psi\,dE_1T\,dE_2R\,dE_3\in\bS_p
\ey
and
\bay
\label{pq}
T\in\bS_p,\!\quad\! R\in\bS_q,\quad \frac1p+\frac1q\le1\!
\!\quad\Longrightarrow\!\quad 
\iiint\!\Psi dE_1TdE_2RdE_3\in\bS_r,\quad\!\frac1r=\frac1p+\frac1q.
\ey

Let us also mention that multiple operator integrals were defined in \cite{Pe4} for 
functions $\Psi$ that belong to the so-called {\it integral projective tensor product} of the corresponding $L^\be$ spaces. We refer the reader to \cite{Pe4} for more detail.

We proceed now to the approach to multiple operator integrals based on the Haagerup tensor product of $L^\be$ spaces. We refer the reader to the book \cite{Pi} for detailed \lb information about Haagerup tensor products.
We define the {\it Haagerup tensor product} \lb
$L^\be(E_1)\!\otimes_{\rm h}\!L^\be(E_2)\!\otimes_{\rm h}\!L^\be(E_3)$ as the space of function $\Psi$ of the form
\bay
\label{htr}
\Psi(x_1,x_2,x_3)=\sum_{j,k\ge0}\a_j(x_1)\b_{jk}(x_2)\g_k(x_3),
\ey
where $\a_j$, $\b_{jk}$, and $\g_k$ are measurable functions such that
\bay
\label{ogr}
\{\a_j\}_{j\ge0}\in L_{E_1}^\be(\ell^2), \quad 
\{\b_{jk}\}_{j,k\ge0}\in L_{E_2}^\be({\mathcal B}),\quad\mbox{and}\quad
\{\g_k\}_{k\ge0}\in L_{E_3}^\be(\ell^2),
\ey
where ${\mathcal B}$ is the space of matrices that induce bounded linear operators on $\ell^2$ and this space is equipped with the operator norm. In other words,
$$
\|\{\a_j\}_{j\ge0}\|_{L^\be(\ell^2)}\df
E_1\mbox{-}\ess\sup\left(\sum_{j\ge0}|\a_j(x_1)|^2\right)^{1/2}<\be,
$$
$$
\|\{\b_{jk}\}_{j,k\ge0}\|_{L^\be({\mathcal B})}\df
E_2\mbox{-}\ess\sup\|\{\b_{jk}(x_2)\}_{j,k\ge0}\|_{{\mathcal B}}<\be,
$$
and
$$
\|\{\g_k\}_{k\ge0}\|_{L^\be(\ell^2)}\df
E_3\mbox{-}\ess\sup\left(\sum_{k\ge0}|\g_k(x_3)|^2\right)^{1/2}<\be.
$$
By the sum on the right-hand of \rf{htr} we mean 
$$
\lim_{M,N\to\be}~\sum_{j=0}^N\sum_{k=0}^M\a_j(x_1)\b_{jk}(x_2)\g_k(x_3).
$$
Clearly, the limit exists.

{\it Throughout the paper by $\sum_{j,k\ge0}$, we mean 
$\lim_{M,N\to\be}~\sum_{j=0}^N\sum^M_{k=0}$}.

The norm of $\Psi$ in 
$L^\be\!\otimes_{\rm h}\!L^\be\!\otimes_{\rm h}\!L^\be$ is, by definition, 
the infimum of
$$
\|\{\a_j\}_{j\ge0}\|_{L^\be(\ell^2)}\|\{\b_{jk}\}_{j,k\ge0}\|_{L^\be({\mathcal B})}
\|\{\g_k\}_{k\ge0}\|_{L^\be(\ell^2)}
$$
over all representations of $\Psi$ of the form \rf{htr}.

It is well known that $L^\be\hat\otimes L^\be\hat\otimes L^\be\subset
L^\be\!\otimes_{\rm h}\!L^\be\!\otimes_{\rm h}\!L^\be$. Indeed,
suppose that $\Psi$ is given by \rf{pred} and \rf{norma} holds. Without loss of generality we may assume 
$$
c_n\df\|\f_n\|_{L^\be}\|\psi_n\|_{L^\be}\|\chi_n\|_{L^\be}\ne0
\quad\mbox{for every}\quad n.
$$
We define $\a_j$, $\b_{j,k}$ and $\g_k$ by
$$
\a_j(x_1)=\frac{\sqrt{c_j}}{\|\f_j\|_{L^\be}}\f_j(x_1),\quad
\g_k(x_3)=\frac{\sqrt{c_k}}{\|\chi_k\|_{L^\be}}\chi_j(x_3)
$$
and
$$
\b_{jk}(x_2)=\left\{\begin{array}{ll}\psi_j(x_2)\|\psi_j\|^{-1}_{L^\be},&j=k\\
[.2cm]0,&j\ne k.
\end{array}\right.
$$
Clearly, \rf{htr} holds,
$$
\|\{\a_j\}\|_{L^\be(\ell^2)}
\le\left(\sum_jc_j\right)^{1/2}<\be,\quad
\|\{\g_k\}\|_{L^\be(\ell^2)}
\le\left(\sum_kc_k\right)^{1/2}<\be
$$
and
$$
\big\|\{\b_{jk}(x_2)\}_{j,k\ge0}\big\|_{{\mB}}\le1.
$$

In \cite{JTT} multiple operator integrals were defined for functions in the Haagerup tensor product of $L^\be$ spaces. Let 
$\Psi\in L^\be\!\otimes_{\rm h}\!L^\be\!\otimes_{\rm h}\!L^\be$ and  
suppose that \rf{htr} and \rf{ogr} hold. The triple operator integral \rf{troi} is defined by
\begin{align}
\label{htraz}
\iiint\Psi(x_1,x_2,x_3)&\,dE_1(x_1)T\,dE_2(x_2)R\,dE_3(x_3)\nonumber\\[.2cm]
=&
\sum_{j,k\ge0}\left(\int\a_j\,dE_1\right)T\left(\int\b_{jk}\,dE_2\right)
R\left(\int\g_k\,dE_3\right)\nonumber\\[.2cm]
=&\lim_{M,N\to\be}~\sum_{j=0}^N\sum_{k=0}^M
\left(\int\a_j\,dE_1\right)T\left(\int\b_{jk}\,dE_2\right)
R\left(\int\g_k\,dE_3\right).
\end{align}

For completeness, we give a proof of the following facts:

\begin{thm}
\label{oprtroi}
{\em(i)} The series in {\em\rf{htraz}} converges in the weak operator topology;
\newline 
{\em(ii)} the sum of the series does not depend on the choice of a representation
{\em\rf{htr}};
\newline
{\em(iii)} the following inequality holds:
\bay
\label{opno}
\!\!\!\!\left\|\iiint\Psi(x_1,x_2,x_3)\,dE_1(x_1)T\,dE_2(x_2)R\,dE_3(x_3)\right\|
\le\|\Psi\|_{L^\be\!\otimes_{\rm h}\!L^\be\!\otimes_{\rm h}\!L^\be}
\|T\|\cdot\|R\|.
\ey
\end{thm}

\Pf Consider the spectral measure $E_2$. It is defined on a $\s$-algebra $\Sigma$ of subsets of
$\X_2$.
We can represent our Hilbert space $\h$ as the direct integral
\bay
\label{din}
\h=\int\limits_{\X_2}\bigoplus\sG(x)\,d\mu(x),
\ey
associated with $E_2$.
Here $\mu$ is a finite measure on $\X_2$, $x\mapsto\sG(x)$, is a measurable Hilbert family. The Hilbert space $\h$ consists of measurable functions $f$ such that $f(x)\in\sG(x)$, $x\in\X_2$, and
$$
\|f\|_\h\df\left(\;\int\limits_{\X_2}\|f(x)\|_{\sG(x)}^2\,d\mu(x)\right)^{1/2}<\be.
$$
Finally, for $\D\in\Sigma$, $E(\D)$ is multiplication by the characteristic function of $\D$. We refer the reader to \cite{BS0}, Ch. 7 for an introduction to direct integrals of Hilbert spaces.

Let us show that the series on the right of \rf{htraz} converges in the weak operator topology.
Let $f$ and $g$ be vectors in $\h$. Put
\bay
\label{ukvj}
u_k\df R\left(\int\g_k\,dE_3\right)f\quad\mbox{and}\quad
v_j\df T^*\left(\int\ov{\a_j}\,dE_1\right)g.
\ey
We consider the vectors $v_j$ and $u_k$ as elements of the direct integral \rf{din}, i.e., vector functions on $\X_2$.

We have
\begin{align*}
&\left|\left(\sum_{j,k\ge0}\left(\int\a_j\,dE_1\right)T\left(\int\b_{jk}\,dE_2\right)
R\left(\int\g_k\,dE_3\right)f,g
\right)\right|\\[.2cm]
&=\left|\sum_{j,k\ge0}\left(\left(\int\b_{jk}\,dE_2\right)u_k,v_j\right)\right|
=\left|\sum_{j,k\ge0}\int\limits_{\X_2}\big(\b_{jk}(x)u_k(x),v_j(x)\big)_{\sG(x)}
\,d\mu(x)\right|
\\[.2cm]
&\le\int\limits_{\X_2}
\|\{\b_{jk}(x)\}_{j,k\ge0}\|_{\mB}\cdot\|\{u_k(x)\}_{k\ge0}\|_{\ell^2}
\cdot\|\{v_j(x)\}_{j\ge0}\|_{\ell^2}d\mu(x)
\\[.2cm]
&\le\|\{\b_{jk}\}_{j,k\ge0}\|_{L^\be(\mB)}
\left(~\int\limits_{\X_2}\Big(\sum_{k\ge0}|u_k(x)|^2\Big)d\mu(x)\right)^{1/2}\!
\!\left(~\int\limits_{\X_2}\Big(\sum_{j\ge0}|v_j(x)|^2\Big)d\mu(x)\right)^{1/2}
\\[.2cm]
&=\|\{\b_{jk}\}_{j,k\ge0}\|_{L^\be(\mB)}
\left(\sum_{k\ge0}\|u_k\|^2_\h\right)^{1/2}
\left(\sum_{j\ge0}\|v_j\|^2_\h\right)^{1/2}.
\end{align*}
Keeping \rf{ukvj} in mind, we see that the last expression is equal to
\begin{align*}
&\|\{\b_{jk}\}_{j,k\ge0}\|_{L^\be(\mB)}
\left(\sum_{k\ge0}\left\|R\left(\int\g_k\,dE_3\right)f\right\|^2_\h\right)^{1/2}
\left(\sum_{j\ge0}\left\|T^*\left(\int\ov{\a_j}\,dE_1\right)g\right\|^2_\h\right)^{1/2}
\\[.2cm]
&\le\|\{\b_{jk}\}_{j,k\ge0}\|_{L^\be\!(\mB)}\|R\|\!\cdot\!\|T\|\!
\left(\sum_{k\ge0}\left\|\left(\int\g_k\,dE_3\right)\!f\right\|^2\right)^{1/2}\!\!\!\!
\left(\sum_{j\ge0}\left\|\left(\int\ov{\a_j}\,dE_1\right)\!g\right\|^2\right)^{1/2}\!\!.
\end{align*}
By properties of integrals with respect to spectral measures,
$$
\sum_{k\ge0}\left\|\left(\int\g_k\,dE_3\right)f\right\|^2=
\left(\int\left(\sum_{k\ge0}|\g_k|^2\right)(dE_3f,f)\right)
\le\|\{\g_k\}_{k\ge0}\|^2_{L^\be(\ell^2)}\|f\|^2.
$$
Similarly,
$$
\sum_{j\ge0}\left\|\left(\int\ov{\a_j}\,dE_1\right)g\right\|^2=
\left(\int\left(\sum_{j\ge0}|\a_j|^2\right)(dE_1g,g)\right)
\le\|\{\a_j\}_{j\ge0}\|^2_{L^\be(\ell^2)}\|g\|^2.
$$
This implies that 
\begin{align*}
&\left|\left(\sum_{j,k\ge0}\left(\int\a_j\,dE_1\right)T\left(\int\b_{jk}\,dE_2\right)
R\left(\int\g_k\,dE_3\right)f,g
\right)\right|\\[.2cm]
&\le
\|\{\b_{jk}\}_{j,k\ge0}\|_{L^\be(\mB)}\cdot
\|\{\a_j\}_{k\ge0}\|_{L^\be(\ell^2)}\cdot
\|\{\g_k\}_{k\ge0}\|_{L^\be(\ell^2)}\|f\|\cdot\|g\|.
\end{align*}
It follows that the series \rf{htraz} converges in the weak operator topology and inequality \rf{opno} holds.


Let us show for completeness that sum \rf{htraz} does not depend on the choice of a representation \rf{htr}. Suppose that \rf{ogr} holds and
$$
\sum_{j,k\ge0}\a_j(x_1)\b_{jk}(x_2)\g_k(x_3)=0\quad\mbox{for almost all}\quad
x_1,~x_2,\quad\mbox{and}\quad x_3. 
$$
We have to show that 
\bay
\label{korre}
\sum_{j,k\ge0}\left(\int\a_j\,dE_1\right)T\left(\int\b_{jk}\,dE_2\right)
R\left(\int\g_k\,dE_3\right)=\0.
\ey
Without loss of generality, we may assume that
$$
\sup_{x_1}\|\{\a_j(x_1)\}_{j\ge0}\|_{\ell^2}<\be,\quad
\sup_{x_3}\|\{\g_k(x_3)\}_{k\ge0}\|_{\ell^2}<\be,
$$
and
$$
\sup_{x_2}\|\{\b_{jk}(x_2)\}_{j,k\ge0}\|_{\mB}<\be.
$$
Put
\bay
\label{vk}
\vk_j(x_2,x_3)\df\sum_{k\ge0}\b_{jk}(x_2)\g_k(x_3)
\ey
Clearly, the series on the right of \rf{vk} converges absolutely and uniformly in $x_2$ and $x_3$ and
$$
\sup_{x_2,\,x_3}\|\{\vk_j(x_2,x_3)\}_{j\ge0}\|_{\ell^2}<\be.
$$
We integrate now the identity
$$
\sum_{j\ge0}\a_j(x_1)\vk_j(x_2,x_3)=0
$$
with respect to the spectral measure $E_3$ and get
\begin{align*}
\int\left(\sum_{j\ge0}\a_j(x_1)\vk_j(x_2,x_3)\right)\,dE_3(x_3)&=
\sum_{j\ge0}\a_j(x_1)\int\vk_j(x_2,x_3)\,dE_3(x_3)\\[.2cm]
&=\sum_{j\ge0}\a_j(x_1)\int\sum_{k\ge0}\b_{jk}(x_2)\g_k(x_3)\,dE_3(x_3)\\[.2cm]
&=\sum_{j\ge0}\a_j(x_1)\sum_{k\ge0}\b_{jk}(x_2)\int\g_k(x_3)\,dE_3(x_3)
=\0.
\end{align*}
Let $u$ be a unit vector in our Hilbert space $\h$. We have
\begin{align*}
R&\left(\sum_{j\ge0}\a_j(x_1)\sum_{k\ge0}\b_{jk}(x_2)\int\g_k(x_3)\,dE_3(x_3)\right)u
\\[.2cm]
&=\sum_{j\ge0}\a_j(x_1)\sum_{k\ge0}\b_{jk}(x_2)R\left(\int\g_k(x_3)\,dE_3(x_3)\right)u=\0.
\end{align*}
Putting 
$$
v_k\df R\left(\int\g_k(x_3)\,dE_3(x_3)\right)u,
$$
we find that
$$
\|v_k\|\le\|R\|\cdot\|\g_k\|_{L^\be(E_3)}
$$
and
$$
\sum_{j\ge0}\a_j(x_1)\sum_{k\ge0}\b_{jk}(x_2)v_k=\0\quad\mbox{for almost all}
\quad x_1\quad\mbox{and}\quad x_2.
$$
Put
$$
\o_k(x_1,x_2)\df\sum_{j\ge0}\a_j(x_1)\b_{jk}(x_2)v_k.
$$
It is easy to see that
$$
\sup_{x_1,\,x_2}\,\sum_{k\ge0}\|\o_k(x_1,x_2)\|<\be\quad\mbox{and}\quad
\sum_{k\ge0}\o_k(x_1,x_2)=\0
\quad\mbox{almost everywhere}.
$$
Clearly, for each $x_1$, the function $x_2\mapsto\o_k(x_2)$ belongs to the projective tensor product $L^\be(E_2)\hat\otimes\h$, we can  
integrate the vector-valued function $\o_k$ with respect to the spectral measure
$E_2$ (see subsection 2.4) and obtain
\begin{align*}
\0&=\int \big(dE_2(x_2)\o_k(x_1,x_2)\big)\\[.2cm]
&=
\sum_{j\ge0}\a_j(x_1)\int\left(dE_2(x_2)\sum_{k\ge0}\b_{jk}(x_2)v_k\right)\\[.2cm]
&=\sum_{j\ge0}\a_j(x_1)\sum_{k\ge0}\left(\int\b_{jk}(x_2)\,dE_2(x_2)\right)v_k
\quad\mbox{for almost all}\quad x_1.
\end{align*}
Thus
\begin{align*}
T&\left(\sum_{j\ge0}\a_j(x_1)\sum_{k\ge0}
\left(\int\b_{jk}(x_2)\,dE_2(x_2)\right)v_k\right)\\[.2cm]
&=
\sum_{j\ge0}\a_j(x_1)\,T
\int \left(dE_2(x_2)\Big(\sum_{k\ge0}\b_{jk}(x_2)\Big)v_k\right)=\0.
\end{align*}
Consider the vectors $w_j$ defined by
$$
w_j\df T\int \left(dE_2(x_2)\Big(\sum_{k\ge0}\b_{jk}(x_2)\Big)v_k\right).
$$
It is easy to see that 
$$
\sum_{j\ge0}\|w_j\|^2<\be.
$$
Integrating the equality 
$$
\sum_{j\ge0}\a_j(x_1)w_j=\0
$$
with respect to the spectral measure $E_1$, we obtain
\begin{align*}
\0&=\sum_{j\ge0}\left(\int\a_j(x_1)\,dE_1(x_1)\right)w_j\\[.2cm]
&=\sum_{j,k\ge0}\left(\int\a_j(x_1)\,dE_1(x_1)\right)T
\left(\int\b_{jk}(x_2)\,dE_2(x_2)\right)R
\left(\int\g_k(x_3)\,dE_3(x_3)\right)u
\end{align*}
which proves \rf{korre}. $\bl$

Note that if $\Psi$ belongs to the projective tensor product 
$L^\be(E_1)\hat\otimes L^\be(E_2)\hat\otimes L^\be(E_3)$, then
the two definitions given above lead to the same result.

It turns out, however, that unlike in the case when the integrand belongs to the projective tensor product $L^\be\hat\otimes L^\be\hat\otimes L^\be$, triple operator integrals with integrands in the Haagerup tensor product 
$L^\be\!\otimes_{\rm h}\!L^\be\!\otimes_{\rm h}\!L^\be$ do not possess property
\rf{bep} with $p<2$; this will be established in \S~\ref{2c}. As for property \rf{pq}, we will show in \S~\ref{SvNtoi} that for integrands in 
$L^\be\!\otimes_{\rm h}\!L^\be\!\otimes_{\rm h}\!L^\be$, property \rf{pq} holds under the assumption $1/p+1/q\le1/2$. We do not know whether \rf{pq} can hold 
if $1/p+1/q>1/2$.

\

\section{\bf Schatten--von Neumann properties of triple operator integrals}
\setcounter{equation}{0}
\label{SvNtoi}

\

In this section we study Schatten--von Nemann properties of triple operator integrals with integrands in the Haagerup tensor product 
$L^\be\!\otimes_{\rm h}\!L^\be\!\otimes_{\rm h}\!L^\be$.
First, we consider the case when one of the operators is bounded and the other one belongs to the Hilbert--Schmidt class. Then we use an interpolation theorem for bilinear operators to a considerably more general situation.

\begin{thm}
\label{bS2}
Let $E_1$, $E_2$, and $E_3$ be spectral measures on Hilbert space and let $\Phi$
be a function in the Haagerup tensor product 
$L^\be(E_1)\!\otimes_{\rm h}\!L^\be(E_2)\!\otimes_{\rm h}\!L^\be(E_3)$.
Suppose that $T$ is a bounded linear operator and $R$ is an operator that belongs 
to the Hilbert--Schmidt class $\bS_2$. Then
\bay
\label{W}
W\df\int\limits_{\X_1}\int\limits_{\X_2}\int\limits_{\X_3}
\Psi(x_1,x_2,x_3)\,dE_1(x_1)T\,dE_2(x_2)R\,dE_3(x_3)\in\bS_2
\ey
and 
\bay
\label{WS_2}
\|W\|_{\bS_2}\le
\|\Psi\|_{L^\be\!\otimes_{\rm h}\!L^\be\!\otimes_{\rm h}\!L^\be}
\|T\|\cdot\|R\|_{\bS_2}.
\ey
\end{thm}

It is easy to see that Theorem \ref{bS2} implies the following fact:

\begin{cor}
\label{S2b}
Let $E_1$, $E_2$, $E_3$, and $\Psi$ satisfy the hypotheses of Theorem 
{\em\ref{bS2}}. If $T$ is a Hilbert Schmidt operator and $R$ is a bounded linear operator, then the operator $W$ defined by {\em\rf{W}} belongs to $\bS_2$ and
$$
\|W\|_{\bS_2}\le
\|\Psi\|_{L^\be\!\otimes_{\rm h}\!L^\be\!\otimes_{\rm h}\!L^\be}
\|T\|_{\bS_2}\|R\|.
$$
\end{cor}

Clearly, to deduce Corollary \ref{S2b} from Theorem \ref{bS2}, it suffices to consider the adjoint operator $W^*$.

\medskip

{\bf Proof of Theorem \ref{bS2}.} Consider first the case when $E_3$ is a discrete spectral measure. In other words, there exists an orthonormal basis 
$\{e_m\}_{m\ge0}$, the spectral measure $E_3$ is defined on the $\s$-algebra of all subsets
of $\Z_+$, and 
$E_3(\{m\})$ is the orthogonal projection onto the one-dimensional space spanned by $e_m$.
In this case the function $\Psi$ has the form
$$
\Psi(x_1,x_2,m)=\sum_{j,k\ge0}\a_j(x_1)\b_{jk}(x_2)\g_k(m),\quad
x_1\in\X_1,~x_2\in\X_2,~m\in\Z_+,
$$
where 
$$
\{\a_j\}_{j\ge0}\in L_{E_1}^\be(\ell^2), \quad 
\{\b_{jk}\}_{j,k\ge0}\in L_{E_2}^\be({\mathcal B}),
$$
and
$$
\sup_{m\ge1}\sum_{k\ge0}|\g_k(m)|^2<\be.
$$
Then
$$
W=\iint\sum_{m\ge0}\Psi(x_1,x_2,m)\,dE_1(x_1)T\,dE_2(x_2)R\,(\cdot,e_m)e_m.
$$
We have
\bay
\label{Zn}
\|W\|_{\bS_2}^2=\sum_{m\ge0}\|We_m\|^2=\sum_{m\ge0}\|Z_mRe_m\|^2,
\ey
where
\begin{align*}
Z_m&\df\iint\Psi(x_1,x_2,m)\,dE_1(x_1)T\,dE_2(x_2)\\[.2cm]
&=\iiint\Psi_m(x_1,x_2,m)\,dE_1(x_1)T\,dE_2(x_2)I\,d\E_m.
\end{align*}
Here $\E_m$ is the spectral measure defined on the one point set $\{m\}$
and the function $\Psi_m$ is defined on $\X_1\times\X_2\times\{m\}$ by
$$
\Psi_m(x_1,x_2,m)=\Psi(x_1,x_2,m),\quad x_1\in\X_1,\quad x_2\in\X_2.
$$

It is easy to see that 
$$
\|\Psi_m\|_{L^\be(E_1)\otimes_{\rm h}L^\be(E_2)\otimes_{\rm h}L^\be(\E_m)}
\le\|\Psi\|_{L^\be(E_1)\otimes_{\rm h}L^\be(E_2)\otimes_{\rm h}L^\be(E_3)},
\quad m\ge0.
$$
It follows now from \rf{opno} that 
$$
\|Z_m\|\le\|\Psi\|_{L^\be\otimes_{\rm h}L^\be\otimes_{\rm h}L^\be}\|T\|,
$$
and by \rf{Zn}, we obtain
\begin{align*}
\sum_{m\ge0}\|We_m\|^2&\le\sum_{n\ge0}\|Z_m\|^2\|Re_m\|^2\\[.2cm]
&\le\|\Psi\|^2_{L^\be\otimes_{\rm h}L^\be\otimes_{\rm h}L^\be}
\|T\|^2\sum_{m\ge0}\|Re_m\|^2
\\[.2cm]
&=\|\Psi\|^2_{L^\be\otimes_{\rm h}L^\be\otimes_{\rm h}L^\be}\|T\|^2\|R\|_{\bS_2}^2.
\end{align*}
It follows that $W\in\bS_2$ and inequality \rf{WS_2} holds.

Consider now the general case. For $N\ge1$, we define the function $\Psi_{[N]}$ by
$$
\Psi_{[N]}\df\sum_{j=0}^N\sum_{k=0}^N
\a_j(x_1)\b_{jk}(x_2)\g_k(x_3).
$$
Since the series on the right-hand side of \rf{htraz} converges weakly, it suffices
to prove that the operators 
$$
W_N\df\iiint
\Psi_{[N]}(x_1,x_2,x_3)\,dE_1(x_1)T\,dE_2(x_2)R\,dE_3(x_3)
$$
belong to $\bS_2$ and 
$$
\|W_N\|_{\bS_2}\le\|\Psi_{[N]}\|_{L^\be\otimes_{\rm h}L^\be\otimes_{\rm h}L^\be}
\|T\|\cdot\|R\|_{\bS_2}
$$
because, obviously,
$$
\|\Psi_{[N]}\|_{L^\be\otimes_{\rm h}L^\be\otimes_{\rm h}L^\be}\le
\|\Psi\|_{L^\be\otimes_{\rm h}L^\be\otimes_{\rm h}L^\be}.
$$
In other words, in the representation of $\Psi$ in the form \rf{htr} we may assume that the sum is finite. We have
$$
W=\sum_{j,k}
\left(\int\a_j\,dE_1\right)T\left(\int\b_{jk}\,dE_2\right)R\left(\int\g_k\,dE_3\right).
$$

We can approximate the functions $\g_k$ by sequences $\g_k^{[n]}$ such that each function $\g_k^{[n]}$ takes at most countably many values, 
$$
\big|\g_k^{[n]}(x)\big|\le|\g_k(x)|,\quad x\in\X_3,
$$
and
$$
\lim_{n\to\be}\big\|\g_k^{[n]}-\g_k\big\|_{L^\be(E_3)}=0.
$$
Consider the operator
$$
W^{[n]}\df\sum_{j,k}
\left(\int\a_j\,dE_1\right)T\left(\int\b_{jk}\,dE_2\right)R\left(\int\g_k^{[n]}\,dE_3\right).
$$
Clearly, in the above representation of $W^{[n]}$ we can replace the spectral measure $E_3$ with a discrete spectral measure whose atoms are the sets on which the 
functions $\g_k^{[n]}$ are constant.

Since we have already proved the desired result in the case when $E_3$ is a discrete spectral measure, we can conclude that $W^{[n]}\in\bS_2$ and
$$
\big\|W^{[n]}\big\|_{\bS_2}\le
\|\Psi\|_{L^\be\otimes_{\rm h}L^\be\otimes_{\rm h}L^\be}\|T\|\cdot\|R\|_{\bS_2}.
$$
To complete the proof, it suffices to observe that 
$$
\lim_{n\to\be}\int\g_k^{[n]}\,dE_3=\int\g_k\,dE_3
$$
in the operator norm. $\bl$

We are going to use Theorem 4.4.1 from \cite{BL} on complex interpolation of bilinear operators. Recall that the Schatten--von Neumann classes $\bS_p$, $p\ge1$, and the space of bounded linear operators ${\mathcal B}(\h)$ form a complex interpolation scale: 
\bay
\label{Spint}
(\bS_1,{\mathcal B}(\h))_{[\theta]}=\bS_{\frac1{1-\theta}},\quad0<\theta<1.
\ey 
This fact is well known. For example, it follows from Theorem 13.1 of Chapter III of \cite{GK}.

\begin{thm}
\label{SNSp}
Let $\Psi\in L^\be(E_1)\!\otimes_{\rm h}\!L^\be(E_2)\!\otimes_{\rm h}\!L^\be(E_3)$.
Then the following holds:
\newline
{\em(i)} if $p\ge2$, $T\in\B(\h)$, and $R\in\bS_p$, then the triple operator integral in {\em\rf{W}} belongs to $\bS_p$ and
\bay
\label{boSp}
\|W\|_{\bS_p}\le\|\Psi\|_{L^\be\!\otimes_{\rm h}\!L^\be\!\otimes_{\rm h}\!L^\be}
\|T\|\cdot\|R\|_{\bS_p};
\ey
\newline
{\em(ii)} if $p\ge2$, $T\in\bS_p$, and $R\in\B(\h)$, then the triple operator integral in {\em\rf{W}} belongs to $\bS_p$ and
$$
\|W\|_{\bS_p}\le\|\Psi\|_{L^\be\!\otimes_{\rm h}\!L^\be\!\otimes_{\rm h}\!L^\be}
\|T\|_{\bS_p}\|R\|;
$$
\newline
{\em(iii)} if $1/p+1/q\le1/2$, $T\in\bS_p$, and $R\in\bS_q$, then the triple operator integral in {\em\rf{W}} belongs to $\bS_r$
with $1/r=1/p+1/q$ and
$$
\|W\|_{\bS_r}\le\|\Psi\|_{L^\be\!\otimes_{\rm h}\!L^\be\!\otimes_{\rm h}\!L^\be}
\|T\|_{\bS_p}\|R\|_{\bS_q}.
$$
\end{thm}

\medskip

We will prove in \S~\ref{2c} that neither (i) nor (ii) holds for $p<2$.

\medskip

{\bf Proof of Theorem \ref{SNSp}.} Let us first prove (i). Clearly, to deduce (ii) from (i), it suffices to consider $W^*$. 

Consider the bilinear operator $\W$ defined by
$$
\W(T,R)=\iiint\Psi(x_1,x_2,x_3)\,dE_1(x_1)T\,dE_2(x_2)R\,dE_3(x_3).
$$
By \rf{opno}, $\W$ maps $\mB(\h)\times\mB(\h)$ into $\mB(\h)$ and
$$
\|\W(T,R)\|\le\|T\|\cdot\|R\|.
$$
On the other hand, by Theorem \ref{bS2}, $\W$ maps $\mB(\h)\times\bS_2$
into $\bS_2$ and
$$
\|\W(T,R)\|_{\bS_2}\le\|T\|\cdot\|R\|_{\bS_2}.
$$
It follows from the complex interpolation theorem for linear operators 
(see \cite{BL}, Theorem 4.1.2 that) $\W$ maps $\mB(\h)\times\bS_p$, $p\ge2$,
into $\bS_p$ and
$$
\|\W(T,R)\|_{\bS_p}\le\|T\|\cdot\|R\|_{\bS_p}.
$$

Suppose now that $1/p+1/q\le1/2$ and $1/r=1/p+1/q$. It follows from statements (i) and (ii) (which we have already proved) that $\W$ maps $\mB(\h)\times\bS_r$ into
$\bS_r$ and $\bS_r\times\mB(\h)$ into $\bS_r$, and
$$
\|\W(T,R)\|_{\bS_r}\le\|T\|\cdot\|R\|_{\bS_r}
\quad\mbox{and}\quad
\|\W(T,R)\|_{\bS_r}\le\|T\|_{\bS_r}\cdot\|R\|.
$$
It follows from  Theorem 4.4.1 of \cite{BL} on interpolation of bilinear operators,
$\W$ maps $(\mB(\h),\bS_r)_{[\theta]}\times(\bS_r,\mB(\h))_{[\theta]}$ into
$\bS_r$ and 
$$
\|\W(T,R)\|_{\bS_r}\le
\|T\|_{(\mB(\h),\bS_r)_{[\theta]}}\|R\|_{(\bS_r,\mB(\h))_{[\theta]}}.
$$
It remains to observe that for $\theta=r/p$,
$$
(\mB(\h),\bS_r)_{[\theta]}=\bS_p\quad\mbox{and}\quad
(\bS_r,\mB(\h))_{[\theta]}=\bS_q,
$$
which is a consequence of \rf{Spint}. $\bl$

\

\section{\bf Haagerup-like tensor products and triple operator integrals}
\setcounter{equation}{0}
\label{Ttoi}

\

We are going to obtain Lipschitz type estimates in the norm of $\bS_p$, $1\le p\le2$, for functions of noncommuting self-adjoint operators in \S~\ref{ple2}. As we have mentioned in the introduction, we are going to use a representation of $f(A_1,B_1)-f(A_2,B_2)$
in terms of triple operator integrals that involve the divided differences
$\dg^{[1]}f$ and $\dg^{[2]}f$. However, we will see in \S~\ref{2c} that the divided differences $\dg^{[1]}f$ and $\dg^{[2]}f$ do not have to belong to the
Haagerup tensor product $L^\be\!\otimes_{\rm h}\!L^\be\!\otimes_{\rm h}\!L^\be$ for an arbitrary function $f$ in the Besov class $B_{\be,1}^1(\R^2)$. In addition to this, representation \rf{predst} involve operators of class $\bS_p$ with $p\le2$. However, we will see in \S~\ref{2c} that statements (i) and (ii) of Theorem \ref{SNSp} do not hold for $p<2$.

To overcome these problems, we offer a new approach to triple operator integrals. In this section we introduce Haagerup-like tensor products and define triple operator integrals whose integrands belong to such Haagerup-like tensor products.

\medskip

{\bf Definition 1.} 
{\it A function $\Psi$ is said to belong to the Haagerup-like tensor product 
$L^\be(E_1)\!\otimes_{\rm h}\!L^\be(E_2)\!\otimes^{\rm h}\!L^\be(E_3)$ of the first kind if it admits a representation
\bay
\label{yaH}
\Psi(x_1,x_2,x_3)=\sum_{j,k\ge0}\a_j(x_1)\b_{k}(x_2)\g_{jk}(x_3),\quad x_j\in\X_j,
\ey
with $\{\a_j\}_{j\ge0},~\{\b_k\}_{k\ge0}\in L^\be(\ell^2)$ and 
$\{\g_{jk}\}_{j,k\ge0}\in L^\be(\mB)$. As usual, 
$$
\|\Psi\|_{L^\be\otimes_{\rm h}\!L^\be\otimes^{\rm h}\!L^\be}
\df\inf\big\|\{\a_j\}_{j\ge0}\big\|_{L^\be(\ell^2)}
\big\|\{\b_k\}_{k\ge0}\big\|_{L^\be(\ell^2)}
\big\|\{\g_{jk}\}_{j,k\ge0}\big\|_{L^\be(\mB)},
$$
the infimum being taken over all representations of the form {\em\rf{yaH}}}.

\medskip

Let us now define triple operator integrals whose integrand belong to the tensor product
$L^\be(E_1)\!\otimes_{\rm h}\!L^\be(E_2)\!\otimes^{\rm h}\!L^\be(E_3)$.

Let $1\le p\le2$. For 
$\Psi\in L^\be(E_1)\!\otimes_{\rm h}\!L^\be(E_2)\!\otimes^{\rm h}\!L^\be(E_3)$, for a bounded linear operator $R$, and for ans operator $T$ of class $\bS_p$, we define the triple operator integral
\bay
\label{WHft}
W=\iint\!\!\upint\Psi(x_1,x_2,x_3)\,dE_1(x_1)T\,dE_2(x_2)R\,dE_3(x_3)
\ey
as the following continuous linear functional on $\bS_{p'}$,
$1/p+1/p'=1$ (on the class of compact operators in the case $p=1$):
\bay
\label{fko}
Q\mapsto
\trace\left(\left(
\iiint
\Psi(x_1,x_2,x_3)\,dE_2(x_2)R\,dE_3(x_3)Q\,dE_1(x_1)
\right)T\right).
\ey

\medskip

Clearly, the triple operator integral in \rf{fko} is well defined because the function
$$
(x_2,x_3,x_1)\mapsto\Psi(x_1,x_2,x_3)
$$ 
belongs to the Haagerup tensor product 
$L^\be(E_2)\!\otimes_{\rm h}\!L^\be(E_3)\!\otimes_{\rm h}\!L^\be(E_1)$. It follows easily from statement (i) of Theorem \ref{SNSp} that
$$
\|W\|_{\bS_p}\le\|\Psi\|_{L^\be\otimes_{\rm h}\!L^\be\otimes^{\rm h}\!L^\be}
\|T\|_{\bS_p}\|R\|,\quad1\le p\le2,
$$
(see Theorem \ref{ftHtp}).

It is easy to see that in the case when $\Psi$ belongs to the projective tensor product $L^\be(E_1)\hat\otimes L^\be(E_2)\hat\otimes L^\be(E_3)$, the definition of the triple operator integral given above is consistent with the definition of the triple operator integral given in \rf{otoi}. Indeed, it suffices to verify this for functions $\Psi$ of the form
$$
\Psi(x_1,x_2,x_3)=\f(x_1)\psi(x_2)\chi(x_3),\quad\f\in L^\be(E_1),\quad
\psi\in L^\be(E_2),\quad\chi\in L^\be(E_3),
$$
in which case the verification is obvious.

We also need trace class triple operator integrals in the case when $T$ is a bounded linear operator and $R\in\bS_p$, $1\le p\le2$.

\medskip

{\bf Definition 2.} {\it A function is said to belong to the Haagerup-like tensor product $L^\be(E_1)\!\otimes^{\rm h}\!L^\be(E_2)\!\otimes_{\rm h}\!L^\be(E_3)$
of the second kind if
$\Psi$ admits a representation
\bay
\label{preds}
\Psi(x_1,x_2,x_3)=\sum_{j,k\ge0}\a_{jk}(x_1)\b_{j}(x_2)\g_k(x_3)
\ey
where $\{\b_j\}_{j\ge0},~\{\g_k\}_{k\ge0}\in L^\be(\ell^2)$, 
$\{\a_{jk}\}_{j,k\ge0}\in L^\be(\mB)$. The norm of $\Psi$ in 
the space $L^\be\otimes^{\rm h}\!L^\be\otimes_{\rm h}\!L^\be$ is defined by
$$
\|\Psi\|_{L^\be\otimes^{\rm h}\!L^\be\otimes_{\rm h}\!L^\be}
\df\inf\big\|\{\a_j\}_{j\ge0}\big\|_{L^\be(\ell^2)}
\big\|\{\b_k\}_{k\ge0}\big\|_{L^\be(\ell^2)}
\big\|\{\g_{jk}\}_{j,k\ge0}\big\|_{L^\be(\mB)},
$$
the infimum being taken over all representations of the form {\em\rf{preds}}}.

\medskip

Suppose now that 
$\Psi\in L^\be(E_1)\!\otimes^{\rm h}\!L^\be(E_2)\!\otimes_{\rm h}\!L^\be(E_3)$,
$T$ is a bounded linear operator, and $R\in\bS_p$, $1\le p\le2$. The continuous linear functional 
$$
Q\mapsto
\trace\left(\left(
\iiint\Psi(x_1,x_2,x_3)\,dE_3(x_3)Q\,dE_1(x_1)T\,dE_2(x_2)
\right)R\right)
$$
on the class $\bS_{p'}$ (on the of compact operators in the case $p=1$) 
determines an operator $W$ of class $\bS_p$, which
we call the triple operator integral
\bay
\label{WHst}
W=\upint\!\!\!\iint\Psi(x_1,x_2,x_3)\,dE_1(x_1)T\,dE_2(x_2)R\,dE_3(x_3).
\ey

Moreover,
$$
\|W\|_{\bS_p}\le
\|\Psi\|_{L^\be\otimes^{\rm h}\!L^\be\otimes_{\rm h}\!L^\be}
\|T\|\cdot\|R\|_{\bS_p}.
$$

As above, in the case when 
$\Psi\in L^\be(E_1)\hat\otimes L^\be(E_2)\hat\otimes L^\be(E_3)$, the definition of the triple operator integral given above is consistent with the definition of the triple operator integral given in \rf{otoi}.

We deduce from Theorem \ref{SNSp} the following Schatten--von Nemann properties of
the triple operator integrals introduced above.

\begin{thm}
\label{ftHtp}
Let $\Psi\in L^\be\!\otimes_{\rm h}\!L^\be\!\otimes^{\rm h}\!L^\be$.
Suppose that $T\in\bS_p$ and $R\in\bS_q$, where
$1\le p\le2$ and $1/p+1/q\le1$. Then the operator $W$ in {\em\rf{WHft}} belongs to $\bS_r$, $1/r=1/p+1/q$, and
\bay
\label{rpq}
\|W\|_{\bS_r}\le\|\Psi\|_{L^\be\otimes_{\rm h}\!L^\be\otimes^{\rm h}\!L^\be}
\|T\|_{\bS_p}\|R\|_{\bS_q}.
\ey
If $T\in\bS_p$, $1\le p\le2$, and $R$ is a bounded linear operator, then $W\in\bS_p$
and
\bay
\label{pB}
\|W\|_{\bS_p}\le\|\Psi\|_{L^\be\otimes_{\rm h}\!L^\be\otimes^{\rm h}\!L^\be}
\|T\|_{\bS_p}\|R\|.
\ey
\end{thm}

\Pf Let $\Phi$ be the function defined by
$$
\Phi(x_2,x_3,x_1)=\Psi(x_1,x_2,x_3).
$$

Consider the case when $R\in\bS_q$, $q\ge1$.
Clearly, the norm of $W$ if $\bS_r$ is the norm of the linear functional 
\rf{fko} on $\bS_{r'}$ (on the class of compact operators if $r=1$). We have
$$
\left|\trace\left(\left(
\iiint
\Psi\,dE_2R\,dE_3Q\,dE_1
\right)T\right)\right|\le\|T\|_{\bS_p}
\left\|\iiint\Psi\,dE_2R\,dE_3Q\,dE_1\right\|_{\bS_{p'}}
$$
(in the case when $p=1$ we have to replace the norm in $\bS_{p'}$ on the right-hand side of the inequality with the operator norm). By Theorem \ref{SNSp},
\begin{align*}
\left\|\iiint\Psi\,dE_2R\,dE_3Q\,dE_1\right\|_{\bS_{p'}}&=
\left\|\iiint\Phi(x_2,x_3,x_1)\,dE_2(x_2)R\,dE_3(x_2)Q\,dE_1(x_1)\right\|_{\bS_{p'}}
\\[.2cm]
&\le\|\Phi\|_{L^\be\!\otimes_{\rm h}\!L^\be\!\otimes_{\rm h}\!L^\be}
\|R\|_{\bS_q}\|Q\|_{\bS_{r'}}\\[.2cm]
&=\|\Psi\|_{L^\be\!\otimes_{\rm h}\!L^\be\!\otimes^{\rm h}\!L^\be}
\|R\|_{\bS_q}\|Q\|_{\bS_{r'}},
\end{align*}
which implies \rf{rpq}. Again, if $p=1$ the norm in $\bS_{p'}$ has to be replaced with the operator norm.

The proof of \rf{pB} is the same. $\bl$

In the same way we can prove the following theorem:

\begin{thm}
\label{stHtp}
Let $\Psi\in L^\be\!\otimes^{\rm h}\!L^\be\!\otimes_{\rm h}\!L^\be$.
Suppose that $p\ge1$, $1\le q\le2$,  and $1/p+1/q\le1$. If $T\in\bS_p$, $R\in\bS_q$, then the operator $W$ in {\em\rf{WHst}} belongs to $\bS_r$, $1/r=1/p+1/q$, and
$$
\|W\|_{\bS_r}\le\|\Psi\|_{L^\be\otimes_{\rm h}\!L^\be\otimes^{\rm h}\!L^\be}
\|T\|_{\bS_p}\|R\|_{\bS_q}.
$$
If $T$ is a bounded linear operator and $R\in\bS_p$, $1\le p\le2$, then $W\in\bS_p$
and
$$
\|W\|_{\bS_p}\le\|\Psi\|_{L^\be\otimes_{\rm h}\!L^\be\otimes^{\rm h}\!L^\be}
\|T\|_{\bS_p}\|R\|.
$$
\end{thm}

\

\section{\bf When do the divided differences $\bs{\dg^{[1]}f}$ and 
$\bs{\dg^{[2]}}f$\\ belong to Haagerup-like tensor products?}
\setcounter{equation}{0}
\label{ddiff}

\

As we have already mentioned before, for functions $f$ in $B_{\be,1}^1(\R^2)$,
the divided differences $\dg^{[1]}f$ and $\dg^{[2]}f$,
$$
\dg^{[1]}f(x_1,x_2,y)\df\frac{f(x_1,y)-f(x_2,y)}{x_1-x_2}
\quad\mbox{and}\quad
\dg^{[2]}f(x,y_1,y_2)\df\frac{f(x,y_1)-f(x,y_2)}{y_1-y_2},
$$
do not have to belong to the Haagerup tensor product $L^\be\!\otimes_{\rm h}\!L^\be\!\otimes_{\rm h}\!L^\be$. This will be prove in \S~\ref{2c}.

In this section we show that for $f\in B_{\be,1}^1(\R^2)$, the divided difference
$\dg^{[1]}f$
belongs to the tensor product 
$L^\be(E_1)\!\otimes_{\rm h}\!L^\be(E_2)\!\otimes^{\rm h}\!L^\be(E_3)$,
while the divided difference 
$\dg^{[2]}f$
belongs to the tensor product 
$L^\be(E_1)\!\otimes^{\rm h}\!L^\be(E_2)\!\otimes_{\rm h}\!L^\be(E_3)$
for arbitrary Borel spectral measures $E_1$, $E_2$, and $E_3$ on $\R$.

This will allow us to prove in the next section that if $(A_1,B_1)$ and $(A_2,B_2)$ are pairs of self-adjoint operators on Hilbert space, $(A_2,B_2)$ is an $\bS_p$ perturbation of $(A_1,B_1)$, $1\le p\le2$, and $f\in B_{\be,1}^1(\R^2)$, then the following integral formula holds:
\begin{align*}
f(A_1,B_1)-f(A_2,B_2)=&
\iint\!\!\upint\frac{f(x_1,y)-f(x_2,y)}{x_1-x_2}
\,dE_{A_1}(x_1)(A_1-A_2)\,dE_{A_2}(x_2)\,dE_{B_1}(y),\\[.2cm]
+&\upint\!\!\!\iint\frac{f(x,y_1)-f(x,y_2)}{y_1-y_2}
\,dE_{A_2}(x)\,dE_{B_1}(y_1)(B_1-B_2)\,dE_{B_2}(y_2).
\end{align*}
 
The following theorem contains a formula that is crucial for our estimates.

\begin{thm}
\label{rrp}
Let $f$ be a bounded function on $\R^2$ whose Fourier transform is supported 
in the ball $\{\xi\in\R^2:~\|\xi\|\le1\}$. Then
\bay
\label{crfor}
\frac{f(x_1,y)-f(x_2,y)}{x_1-x_2}=
\sum_{j,k\in\Z}\frac{\sin(x_1-j\pi)}{x_1-j\pi}\cdot\frac{\sin(x_2-k\pi)}{x_2-k\pi}
\cdot\frac{f(j\pi,y)-f(k\pi,y)}{j\pi-k\pi},
\ey
where for $j=k$, we assume that 
$$
\frac{f(j\pi,y)-f(k\pi,y)}{j\pi-k\pi}
=\frac{\partial f(x,y)}{\partial x}\Big|_{(j\pi,y)}.
$$

Moreover,
\bay
\label{sinusy}
\sum_{j\in\Z}\frac{\sin^2(x_1-j\pi)}{(x_1-j\pi)^2}
=\sum_{k\in\Z}\frac{\sin^2(x_2-k\pi)}{(x_2-k\pi)^2}=1,
\quad x_1~x_2\in\R,
\ey
and
\bay
\label{comnor}
\sup_{y\in\R}\left\|\left\{\frac{f(j\pi,y)-f(k\pi,y)}{j\pi-k\pi}
\right\}_{j,k\in\Z}\right\|_{\mB}\le\const\|f\|_{L^\be(\R)}.
\ey
\end{thm}

\medskip

To prove the theorem, we are going to use the construction in the proof of  Theorem 6.1 of \cite{APPS}. 

\medskip 

\Pf Given $y\in\R$, we consider the function $f_y$ on $\R$ defined by
$f_y(x)=f(x,y)$. Clearly, $f_y$ is a bounded function whose Fourier transform
is supported in $[-1.1]$.
We apply Theorem 6.1 of \cite{APPS} for $f_y$. By formula (6.4) of \cite{APPS},
we have
\bay
\label{rrr}
\frac{f_y(x_1)-f_y(x_2)}{x_1-x_2}=\
\sum_{k\in\Z}\frac{f_y(x_1)-f_y\big(k\pi\big)}{x_1-k\pi}
\cdot\frac{\sin(x_2-k\pi)}{x_2-k\pi}.
\ey
Moreover, by inequality (6.6) of \cite{APPS},
$$
\sum_{k\in\Z}\frac{|f_y(x_1)-f_y(k\pi)|^2}{(x_1-k\pi)^2}
\le3\|f\|_{L^\be(\R)}^2.
$$
It is well known (see, e.g., \cite{Ti}, 3.3.2, Example IV) that
$$
\sum_{n\in\Z}\frac{\sin^2(x-n\pi)}{(x-n\pi)^2}=1,\quad x\in\R,
$$
and so \rf{sinusy} holds.

It follows that the series on the right-hand side of \rf{rrr} converges pointwise.
Note that on the right-hand side of \rf{rrr} in the case $x_1=k\pi$, we assume that
$$
\frac{f_y(x_1)-f_y\big(k\pi\big)}{x_1-k\pi}=f'_y(k\pi).
$$

Applying formula (6.4) of \cite{APPS} for the second time, we obtain
\bay
\label{fjk}
\frac{f_y(x_1)-f_y\big(k\pi\big)}{x_1-k\pi}=\sum_{j\in\Z}
\frac{f_y(j\pi)-f_y(k\pi)}{j\pi-k\pi}\cdot\frac{\sin(x_1-j\pi)}{x_1-j\pi}.
\ey
Again, in the case $j=k$ we assume that
$$
\frac{f_y(j\pi)-f_y\big(k\pi\big)}{j\pi-k\pi}=f'_y(j\pi)=\frac{\partial f(x,y)}{\partial x}\Big|_{(j\pi,y)}.
$$

Clearly, \rf{crfor} is a consequence of \rf{rrr} and \rf{fjk}.

Let us estimate the operator norm of the matrix
$$
\left\{\frac{f(j\pi,y)-f(k\pi,y)}{j\pi-k\pi}
\right\}_{j,k\in\Z}
$$
We represent this matrix as the sum of the matrices $C_y=\{c_{jk}(y)\}_{j,k\in\Z}$
and $D_y=\{d_{jk}(y)\}_{j,k\in\Z}$, where
$$
c_{jk}(y)=\left\{\begin{array}{ll}\frac{f(j\pi,y)-f(k\pi,y)}{j\pi-k\pi},&j\ne k\\[.2cm]
0,&j=k
\end{array}\right.
$$
and
$$
d_{jk}(y)=\left\{\begin{array}{ll}0,&j\ne k\\[.2cm]
\frac{\partial f(x,y)}{\partial x}\Big|_{(j\pi,y)},&j=k.
\end{array}\right.
$$

To estimate the operator norm of $C_y$, we observe that $C_y$ is the commutator of the discrete Hilbert transform ${\mathcal H}_{\rm d}$ and a multiplication operator on 
$\ell^2$. Recall that the discrete Hilbert transform ${\mathcal H}_{\rm d}$ on the two-sided sequence space $\ell^2_\Z$ is the operator with matrix 
$\{h_{jk}\}_{j,k\in\Z}$ defined by
$$
h_{jk}=\left\{\begin{array}{ll}\frac1{j-k},&j\ne k\\[.2cm]
0,&j=k.
\end{array}\right.
$$
It is well known that ${\mathcal H}_{\rm d}$ is a bounded linear operator
on $\ell^2(\Z)$. Indeed, $h_{jk}=\hat\psi(j-k)$, where $\phi$ is the bounded function on the unit circle $\T$ defined by
$$
\phi(e^{{\rm i}t})={\rm i}(\pi-t),\quad0\le t<2\pi,
$$
(see, e.g., \cite{Pe5}, Ch. I, \S~1). It follows that ${\mathcal H}_{\rm d}$ is bounded because if we identify the two-sided sequence space $\ell^2_\Z$ with the space $L^2(\T)$ via the unitary map
$$
\{c_n\}_{n\in\Z}\mapsto\sum_{n\in\Z}c_nz^n,
$$
the operator ${\mathcal H}_{\rm d}$ becomes the operator of multiplication on $L^2(\T)$, and so it is bounded and $\|{\mathcal H}_{\rm d}\|=\pi$.

It is easy to see that the matrix of $\pi C_y$ coincides with the matrix of the commutator $M_{f_y}{\mathcal H}_{\rm d}-{\mathcal H}_{\rm d}M_{f_y}$ of the discrete Hilbert transform and the multiplication operator $M_{f_y}$ whose matrix is diagonal with diagonal entries $\{f(j\pi,y)\}_{j\in\Z}$. Clearly, 
$$
\|M_{f_y}\|=\sup_{j\in\Z}|f(j\pi,y)|\le\|f\|_{L^\be(\R^2)}.
$$
Thus
$$
\|C_y\|=\frac1\pi\|M_{f_y}{\mathcal H}_{\rm d}-{\mathcal H}_{\rm d}M_{f_y}\|\le
2\|M_{f_y}\|\cdot\|{\mathcal H}_{\rm d}\|\le2\|f\|_{L^\be(\R^2)}.
$$

On the other hand,
$$
\|D_y\|=\sup_{j\in\Z}\left|\frac{\partial f(x,y)}{\partial x}\Big|_{(j\pi,y)}\right|
\le\|f\|_{L^\be(\R^2)}
$$
by Bernstein's inequality. This completes the proof of \rf{comnor}.
$\bl$

\medskip

{\bf Remark.} It is easy to see from the proof of Theorem \ref{rrp} that one can replace the condition $\supp\F f\subset\{\xi\in\R^2:~\|\xi\|\le1\}$ with the condition
$\supp\F f\subset[-1,1]\times\R$.

\begin{cor}
\label{sle}
Let $f$ be a bounded function on $\R^2$ such that its Fourier transform is supported in $\{\xi\in\R^2:~\|\xi\|\le\s\}$, $\s>0$. Then the divided differences
$\dg^{[1]}f$ and $\dg^{[2]}f$ have the following properties: 
$$
\dg^{[1]}f\in L^\be(E_1)\!\otimes_{\rm h}\!L^\be(E_2)\!\otimes^{\rm h}\!L^\be(E_3)
\quad\mbox{and}\quad
\dg^{[2]}f\in L^\be(E_1)\!\otimes^{\rm h}\!L^\be(E_2)\!\otimes_{\rm h}\!L^\be(E_3)
$$
for arbitrary Borel spectral measures $E_1$, $E_2$ and $E_3$. Moreover,
\bay
\label{dg1}
\big\|\dg^{[1]}f\big\|_{L^\be\!\otimes_{\rm h}\!L^\be\!\otimes^{\rm h}\!L^\be}
\le\const\s\|f\|_{L^\be(\R^2)}
\ey
and
\bay
\label{dg2}
\big\|\dg^{[2]}f\big\|_{L^\be\!\otimes^{\rm h}\!L^\be\!\otimes_{\rm h}\!L^\be}
\le\const\s\|f\|_{L^\be(\R^2)}.
\ey
\end{cor}

\Pf Inequality \rf{dg1} in the case $\s=1$ is an immediate consequence of Theorem \ref{rrp}. It is easy to see that by rescaling the function $f$, we obtain inequality \rf{dg1} for an arbitrary positive number $\s$. Inequality \rf{dg2} can be deduced from inequality \rf{dg1}
by applying \rf{dg2} to the function $g$ defined by
$$
g(x_1,x_2,y)=f(y,x_1,x_2).\quad\bl
$$

\begin{thm}
\label{Bes}
Let $f\in B_{\be,1}^1(\R^2)$. Then
$$
\dg^{[1]}f\in L^\be(E_1)\!\otimes_{\rm h}\!L^\be(E_2)\!\otimes^{\rm h}\!L^\be(E_3)
\quad\mbox{and}\quad
\dg^{[2]}f\in L^\be(E_1)\!\otimes^{\rm h}\!L^\be(E_2)\!\otimes_{\rm h}\!L^\be(E_3)
$$
for arbitrary Borel spectral measures $E_1$, $E_2$ and $E_3$. Moreover,
$$
\big\|\dg^{[1]}f\big\|_{L^\be\!\otimes_{\rm h}\!L^\be\!\otimes^{\rm h}\!L^\be}
\le\const\|f\|_{B_{\be,1}^1}
$$
and
$$
\big\|\dg^{[2]}f\big\|_{L^\be\!\otimes^{\rm h}\!L^\be\!\otimes_{\rm h}\!L^\be}
\le\const\s\|f\|_{B_{\be,1}^1}.
$$
\end{thm}

\Pf Let $f_n=f*W_n$, $n\in\Z$ (see Subsection 2.1.1). Then $f_n$ satisfies the hypotheses of Corollary \ref{sle} with $\s=2^{n+1}$.
By Corollary \ref{sle}, we have
\begin{align*}
\big\|\dg^{[1]}f\big\|_{L^\be\!\otimes_{\rm h}\!L^\be\!\otimes^{\rm h}\!L^\be}
&=
\left\|\sum_{n\in\Z}\dg^{[1]}f_n\right\|_{L^\be\!\otimes_{\rm h}\!L^\be\!\otimes^{\rm h}\!L^\be}
\le\sum_{n\in\Z}\big\|\dg^{[1]}f_n\big\|_{L^\be\!\otimes_{\rm h}\!L^\be\!\otimes^{\rm h}\!L^\be}\\[.2cm]
&\le\const\sum_{n\in\Z}2^{n+1}\|f_n\|_{L^\be}\le\const\|f\|_{B_{\be,1}^1}.
\end{align*}
The proof of the result for $\dg f^{[2]}$ is the same. $\bl$

\

\section{\bf Lipschitz type estimates in the case $\bs{1\le p\le2}$}
\setcounter{equation}{0}
\label{ple2}

\

In this section we prove that for functions $f$ in the Besov class 
$B_{\be,1}^1(\R^2)$, we have a Lipschitz type estimate for functions of noncommuting self-adjoint operators in the norm of $\bS_p$ with $p\in[1,2]$. To this end, we first prove the integral formula given in the introduction.

\begin{thm}
\label{intfor}
Let $f\in B_{\be,1}^1(\R^2)$ and $1\le p\le2$. Suppose that $(A_1,B_1)$ and $(A_2,B_2)$ are pairs of self-adjoint operators such that $A_2-A_1\in\bS_p$ and 
$B_2-B_1\in\bS_p$. Then the following identity holds:
\begin{align}
\label{osnf}
f(A_1,B_1)&-f(A_2,B_2)\nonumber\\[.2cm]
&=
\iint\!\!\upint\frac{f(x_1,y)-f(x_2,y)}{x_1-x_2}
\,dE_{A_1}(x_1)(A_1-A_2)\,dE_{A_2}(x_2)\,dE_{B_1}(y),\nonumber\\[.2cm]
&+\upint\!\!\!\iint\frac{f(x,y_1)-f(x,y_2)}{y_1-y_2}
\,dE_{A_2}(x)\,dE_{B_1}(y_1)(B_1-B_2)\,dE_{B_2}(y_2).
\end{align}
\end{thm}

Note that by Theorem \ref{Bes}, the divided differences $\dg^{[1]}f$ and $\dg^{[2]}f$ belong to the corresponding Haagerup like tensor products, and so the triple operator integrals on the right make sense.

\medskip

\Pf It suffices to prove that
\begin{align}
\label{perfor}
f(A_1,B_1)&-f(A_2,B_1)\nonumber\\[.2cm]
&=\iint\!\!\upint\big(\dg^{[1]}f\big)(x_1,x_2,y)
\,dE_{A_1}(x_1)(A_1-A_2)\,dE_{A_2}(x_2)\,dE_{B_1}(y)
\end{align}
and 
\begin{align}
\label{vtfor}
f(A_2,B_1)&-f(A_2,B_2)\nonumber\\[.2cm]
&=\upint\!\!\!\iint\big(\dg^{[2]}f\big)(x,y_1,y_2)
\,dE_{A_2}(x)\,dE_{B_1}(y_1)(B_1-B_2)\,dE_{B_2}(y_2).
\end{align}

Let us establish \rf{perfor}. Formula \rf{vtfor} can be proved in exactly the same way. 

Suppose first that the function $\dg^{[1]}f$ belongs to the projective tensor product \lb$L^\be(E_{A_1})\hat\otimes L^\be(E_{A_2})\hat\otimes L^\be(E_{B_1})$. 
In this case we can write
\begin{align*}
\iint\!\!\upint\big(\dg^{[1]}f\big)(x_1,x_2,y)&
\,dE_{A_1}(x_1)(A_1-A_2)\,dE_{A_2}(x_2)\,dE_{B_1}(y)\\[.2cm]
&=\iint\!\!\upint\big(\dg^{[1]}f\big)(x_1,x_2,y)
\,dE_{A_1}(x_1)A_1\,dE_{A_2}(x_2)\,dE_{B_1}(y)\\[.2cm]
&-\iint\!\!\upint\big(\dg^{[1]}f\big)(x_1,x_2,y)
\,dE_{A_1}(x_1)A_2\,dE_{A_2}(x_2)\,dE_{B_1}(y).
\end{align*}
Note that the above equality does not make sense if 
$\dg^{[1]}f$ does not belong to $L^\be\hat\otimes L^\be\hat\otimes L^\be$ because
the operators $A_1$ and $A_2$ do not have to be compact, while the definition of
triple operator integrals with integrands in the Haagerup-like tensor product
$L^\be\!\otimes_{\rm h}\!L^\be\!\otimes^{\rm h}\!L^\be$
assumes that the operators $A_1$ and $A_2$ belong to $\bS_2$. 

It follows immediately from the definition of triple operator integrals with integrands in $L^\be\hat\otimes L^\be\hat\otimes L^\be$ that
\begin{align*}
\iint\!\!\upint\big(\dg^{[1]}f\big)(x_1,x_2,y)&
\,dE_{A_1}(x_1)A_1\,dE_{A_2}(x_2)\,dE_{B_1}(y)\\[.2cm]
=&\iint\!\!\upint x_1\big(\dg^{[1]}f\big)(x_1,x_2,y)
\,dE_{A_1}(x_1)\,dE_{A_2}(x_2)\,dE_{B_1}(y)
\end{align*}
and
\begin{align*}
\iint\!\!\upint\big(\dg^{[1]}f\big)(x_1,x_2,y)&
\,dE_{A_1}(x_1)A_2\,dE_{A_2}(x_2)\,dE_{B_1}(y)\\[.2cm]
=&\iint\!\!\upint x_2\big(\dg^{[1]}f\big)(x_1,x_2,y)
\,dE_{A_1}(x_1)\,dE_{A_2}(x_2)\,dE_{B_1}(y).
\end{align*}
Thus
\begin{align*}
\iint\!\!\upint\big(\dg^{[1]}f\big)&(x_1,x_2,y)
\,dE_{A_1}(x_1)A_1\,dE_{A_2}(x_2)\,dE_{B_1}(y)\\[.2cm]
-&\iint\!\!\upint\big(\dg^{[1]}f\big)(x_1,x_2,y)
\,dE_{A_1}(x_1)A_2\,dE_{A_2}(x_2)\,dE_{B_1}(y)\\[.2cm]
=&\iint\!\!\upint(x_1-x_2)\frac{f(x_1,y)-f(x_2,y)}{x_1-x_2}
\,dE_{A_1}(x_1)\,dE_{A_2}(x_2)\,dE_{B_1}(y)\\[.2cm]
=&\iint\!\!\upint f(x_1,y)\,dE_{A_1}(x_1)\,dE_{A_2}(x_2)\,dE_{B_1}(y)\\[.2cm]
-&\iint\!\!\upint f(x_2,y)\,dE_{A_1}(x_1)\,dE_{A_2}(x_2)\,dE_{B_1}(y)
=f(A_1,B_1)-f(A_2,B_1).
\end{align*}


As in the proof of Theorem \ref{Bes}, we consider the functions $f_n$ defined by $f_n=f*W_n$,
$n\in\Z$. It is easy to see from the definition and properties of the Besov class $B_{\be,1}^1(\R^2)$ that to prove \rf{perfor}, it suffices to show that
\begin{align*}
f_n(A_1,B_1)&-f_n(A_2,B_1)\\[.2cm]
&=\iint\!\!\upint\big(\dg^{[1]}f_n\big)(x_1,x_2,y)
\,dE_{A_1}(x_1)(A_1-A_2)\,dE_{A_2}(x_2)\,dE_{B_1}(y).
\end{align*}

As we have mentioned in Subsection 2.1.1, $f_n$ is a restriction of an entire function of two variables to 
$\R\times\R$. Thus it suffices to establish formula \rf{perfor} in the case when $f$ is an entire function. To complete the proof, we show that for entire functions $f$ the divided differences $\dg^{[1]}f$ must belong to the projective tensor product
$L^\be(E_{A_1})\hat\otimes L^\be(E_{A_2})\hat\otimes L^\be(E_{B_1})$.  

Let $f(x,y)=\sum\limits_{n=0}^\be(\sum\limits_{m=0}^\be a_{mn}x^my^n)$ be an entire function and let $R$ be a positive number such that the spectra $\s(A_1)$, 
$\s(A_2)$, and $\s(B)$ are contained in $[-R/2,R/2]$.
Clearly,
$$
\|f\|_{L^\be\hat\otimes L^\be}\le\sum\limits_{n=0}^\be\Big(\sum\limits_{m=0}^\be |a_{mn}|R^{m+n}\Big)
<\be
$$
and
\begin{align*}
\left\|\dg^{[1]}f\right\|_{L^\be\hat\otimes L^\be\hat\otimes L^\be}&=
\left\|\sum\limits_{n=0}^\be\left(\sum\limits_{m=1}^\be
\Big(\sum_{j=0}^{m-1} a_{mn}x_1^jx_2^{m-1-j}y^n\Big)
\right)\right\|_{L^\be\hat\otimes L^\be\hat\otimes L^\be}\\[.2cm]
&\le\sum\limits_{n=0}^\be\left(\sum\limits_{m=1}^\be m|a_{mn}|R^{m+n-1}\right)<+\be,
\end{align*}
where in the above expressions $L^\be$ means $L^\be[-R,R]$. This completes the proof. $\bl$


\begin{thm}
\label{Lipp>2}
Let $p\in[1,2]$. Then there is a positive number $C$ such that
\bay
\label{Lipnerp2}
\|f(A_1,B_1)-f(A_2,B_2)\|\le C\|f\|_{B_{\be,1}^1}
\max\big\{\|A_1-A_2\|_{\bS_p},\|B_1-B_2\|_{\bS_p}\big\},
\ey
whenever $f\in B_{\be,1}^1(\R^2)$, and $A_1$, $A_2$, $B_1$, and $B_2$ are self-adjoint operators such that
$A_2-A_1\in\bS_p$ and 
$B_2-B_1\in\bS_p$.
\end{thm}

\Pf This is an immediate consequence of Theorem \ref{intfor} and Theorems
\ref{ftHtp} and \ref{stHtp}. $\bl$

\medskip

{\bf Remark 1.} We have defined functions $f(A,B)$ for $f$ in 
$B_{\be,1}^1(\R^2)$ only for bounded self-adjoint operators $A$ and $B$. 
However, formula \rf{osnf} allows us to define the difference $f(A_1,B_1)-f(A_2,B_2)$ in the case when $f\in B_{\be,1}^1(\R^2)$ and the 
self-adjoint operators $A_1,\,A_2,\,B_1,\,B_2$ are possibly unbounded once we know that the pair $(A_2,B_2)$ is an $\bS_p$ perturbation of the pair $(A_1,B_1)$, 
$1\le p\le2$.
Moreover, inequality \rf{Lipnerp2} also holds for such operators.

\medskip

{\bf Remark 2.} Let $\fI$ be an operator ideal that is an interpolation ideal between $\bS_1$ and $\bS_2$. Then it follows easily from Theorems
\ref{ftHtp} and \ref{stHtp} that for $f\in B_{\be,1}^1(\R^d)$ and for self-adjoint operators $A_1$, $A_2$, $B_1$, $B_2$ with $A_1-A_2\in\fI$, $B_1-B_2\in\fI$, the following inequality holds:
$$
\|f(A_1,B_1)-f(A_2,B_2)\|_\fI\le\const\|f\|_{B_{\be,1}^1}
\max\big\{\|A_1-A_2\|_\fI,\|B_1-B_2\|_\fI\big\}.
$$

To complete the section, we state a problem. 
\medskip

{\bf Problem.} It is well known that if $f$ is an arbitrary Lipschitz function on 
$\R^2$, then the following inequality holds:
\bay
\label{?}
\|f(A_1,B_1)-f(A_2,B_2)\|_{\bS_2}\le\const\|f\|_{\Lip}
\max\big\{\|A_1-A_2\|_{\bS_2},\|B_1-B_2\|_{\bS_2}\big\}
\ey
for arbitrary pairs $(A_1,B_1)$ and $(A_2,B_2)$ of {\it commuting} self-adjoint operators such that $A_1-A_2\in\bS_2$ and $B_1-B_2\in\bS_2$. As we have mentioned in the introduction, the same is true in the Schatten--von Neumann norm $\bS_p$ with $1<p<\be$ which was proved recently in \cite{KPSS}. We do not know whether inequality \rf{?} holds for pairs of noncommuting self-adjoint operators. Certainly, we have not defined functions $f(A,B)$ for all Lipschitz functions $f$ and all pairs of self-adjoint operators $(A,B)$. However, we can consider pairs of finite rank self-adjoint operators $(A_1,B_1)$ and $(A_2,B_2)$ and ask the question of whether inequality \rf{?} holds for such pairs.

\

\section{\bf No Lipschitz type estimates in the operator norm\\ and in the 
$\bs{\bS_p}$ norm for $\bs{p>2}\,\,$ !}
\setcounter{equation}{0}
\label{Bp>2}

\

The purpose of this section is to show that there is no Lipschitz type inequality of the form \rf{Lipp>2} in the norm of $\bS_p$ with $p>2$ and in the operator norm for
an arbitrary function $f$ in $B_{\be,1}^1(\R^2)$. 

\begin{thm} 
\label{nLte}
{\em(i)}
There is no positive number $M$ such that 
$$
\|f(A_1,B)-f(A_2,B)\|\le M\|f\|_{L^\be(\R^2)}\|A_1-A_2\|
$$ 
for all bounded functions $f$ on $\R^2$ with Fourier transform supported 
in $[-2\pi,2\pi]^2$ and for all finite rank self-adjoint operators $A_1,\,A_2,\,B$.

{\em(ii)}
Let $p>2$.
Then there is no positive number $M$ such that 
$$
\|f(A_1,B)-f(A_2,B)\|_{\bS_p}\le M\|f\|_{L^\be(\R^2)}\|A_1-A_2\|_{\bS_p}
$$ 
for all bounded functions $f$ on $\R^2$ with Fourier transform supported in $[-2\pi,2\pi]^2$ and for all finite rank self-adjoint operators $A_1,\,A_2,\,B$.
\end{thm}

\Pf Let us first prove (ii). Let $\{g_j\}_{1\le j\le N}$ and $\{h_j\}_{1\le j\le N}$
be orthonormal systems in Hilbert space. 
Consider the rank one projections $P_j$ and $Q_j$ defined by
$$
P_j v=(v,g_j)g_j\quad\mbox{and}\quad Q_jv=(v,h_j)h_j,\quad 1\le j\le N.
$$
We define the self-adjoint operators $A_1$, $A_2$, and $B$ by
$$
A_1=\sum_{j=1}^N2jP_j,\quad A_2=\sum_{j=1}^N(2j+1) P_j,\quad\mbox{and}\quad
B=\sum_{k=1}^Nk\, Q_k.
$$
Then $\|A_1-A_2\|_{\bS_p}=N^{\frac1p}$.
Put
$$
\f(x)=\frac{1-\cos2\pi x}{2\pi^2 x^2}.
$$
Clearly, 
$\supp\F\f\subset[-2\pi,2\pi]$, $\f(k)=0$ for all $k\in\Z$ such that $k\ne0$, $\f(0)=1$.
Put $\f_k(x)=\f(x-k)$.
Given a matrix $\{\tau_{jk}\}_{1\le j,k\le N}$, we define 
the function $f$ by
$$
f(x,y)=\sum_{1\le j,k\le N}\tau_{jk}\f_{2j}(x)\f_k(y).
$$ 
It is easy to see that $\f_{2j}(A_1)=P_j$, $\f_{2j}(A_2)=0$, $\f_k(B)=Q_k$ provided
$1\le j,k\le N$, and
$$
\|f\|_{L^\be(\R^2)}\le\const\max_{1\le j,k\le N}|\tau_{jk}|.
$$
Clearly,
$$
f(A_1,B)=
\sum_{1\le j,k\le N}\tau_{jk}P_jQ_k\quad\mbox{and}\quad
f(A_2,B)=\0.
$$
Note that 
$$
(f(A_1,B)h_k,g_j)=\tau_{jk}(h_k,g_j),\quad1\le j,k\le N.
$$  
Clearly, for every unitary matrix
$\{u_{jk}\}_{1\le j,k\le N}$, there exist orthonormal systems $\{g_j\}_{1\le j\le N}$ and $\{h_j\}_{1\le j\le N}$
such that  $(h_k,g_j)=u_{jk}$.
Put 
$$
u_{jk}\df\frac1{\sqrt N}\exp\left(\frac{2\pi{\rm i}jk}N\right),\quad 1\le j,k\le N.
$$
Obviously, $\{u_{jk}\}_{1\le j,k\le N}$ is a unitary matrix.
Hence, 
we may find vectors $\{g_j\}_{j=1}^N$ and $\{h_j\}_{j=1}^N$ such that $(h_k,g_j)=u_{jk}$.
Put $\tau_{jk}=\sqrt N \,\,\ov u_{jk}$. Then
$$
\|f(A_1,B)\|_{\bS_p}=\|\{|u_{jk}|\}_{1\le j,k\le N}\|_{\bS_p}=\|\{|u_{jk}|\}_{1\le j,k\le N}\|_{\bS_2}=\sqrt N
$$
because $\rank\{|u_{jk}|\}_{1\le j,k\le N}=1$.
So for each positive integer $N$ we have constructed a function $f$ and operators $A_1$, $A_2$, $B$ 
such that $|f|\le\const$, $\supp\F f\subset[-2\pi,2\pi]^2$, $\|A_1-A_2\|_{\bS_p}=N^{\frac1p}$ and
$\|f(A_1,B)-f(A_2,B)\|_{\bS_p}=\sqrt N$. It remains to observe that $\lim_{N\to\be}N^{\frac12-\frac1p}=\be$
for $p>2$. 

Exactly the same construction works to prove (i). It suffices to replace in the above construction the $\bS_p$ norm with the operator norm and observe that 
$\|A_1-A_2\|=1$ and $\|f(A_1,B)-f(A_2,B)\|=\sqrt N$. $\bl$

Theorem \ref{nLte} implies that there is no Lipschitz type estimate in the operator norm and in the $S_p$ norm with $p>2$. Note that in the construction given in the proof the norms of $A_1-A_2$ cannot get small. The following result shows that we can easily overcome this problem. 

\begin{thm}
\label{BnLte}
There exist a sequence $\{f_n\}_{n\ge0}$ of functions in $B_{\be,1}^1(\R^2)$ and sequences of self-adjoint finite rank operators $\big\{A_1^{(n)}\big\}_{n\ge0}$, 
$\big\{A_2^{(n)}\big\}_{n\ge0}$, 
and $\big\{B^{(n)}\big\}_{n\ge0}$
such that the norms $\|f_n\|_{B_{\be,1}^1}$ do not depend on $n$,
$$
\lim_{n\to\be}\big\|A_1^{(n)}-A_2^{(n)}\big\|\to0,\quad\mbox{but}\quad
\|f(A_1,B)-f(A_2,B)\|\to\be.
$$
The same is true in the norm of $\bS_p$ for $p>2$.
\end{thm}

\Pf The existence of such sequences can be obtained easily from the construction in the proof of 
Theorem \ref{nLte}. It suffices to make the following observation. Let $f$, $A_1$, $A_2$ and $B$ be as in the proof of Theorem \ref{nLte} and let $\e>0$. Put 
$f_\e(x,y)\df\e f\big(\frac{x}{\e},\frac{y}{\e}\big)$. Then 
$$
\|f_\e\|_{B_{\be 1}^1}=\|f\|_{B_{\be 1}^1}, 
\quad \|f_\e(\e A_1,\e B)-f_\e(\e A_2,\e B)\|=\e N^{1/2},
\quad\mbox{and}\quad \|\e A_1-\e A_2\|=\e.
$$
If $p>2$, then
$$
\|f_\e(\e A_1,\e B)-f_\e(\e A_2,\e B)\|_{\bS_p}=\e N^{1/2}
\quad\mbox{and}\quad \|\e A_1-\e A_2\|_{\bS_p}=\e N^{1/p}.\quad\bl
$$

\medskip

{\bf Remark.} The construction given in the proof of Theorem \ref{nLte} shows that for every
positive number $M$ there exist a function $f$ on $\R^2$ whose Fourier transform is supported in 
$[-2\pi,2\pi]^2$ such that $\|f\|_{L^\be(\R)}\le\const$ and self-adjoint operators of finite rank $A_1$, $A_2$, $B$ such that $\|A_1-A_2\|=1$, but $\|f(A_1,B)-f(A_2,B)\|>M$.
It follows that unlike in the case of commuting self-adjoint operators (see \cite{APPS}), the fact that $f$ is a H\"older function of order $\a\in(0,1)$ on $\R^2$ does not imply the H\"older type estimate 
$$
\|f(A_1,B_1)-f(A_2,B_2)\|\le\const\max\big\{\|A_1-A_2\|^\a,\|B_1-B_2\|^\a\big\}.
$$

\

\section{\bf Two counterexamples}
\setcounter{equation}{0}
\label{2c}

\

We apply the results of the previous section to show that statements (i) and (ii) of Theorem \ref{SNSp} do not hold for 
$p\in[1,2)$. We also deduce from the results of \S~\ref{Bp>2} that the divided differences $\dg^{[1]}f$ and $\dg^{[2]}f$ do not have to belong to
the Haagerup tensor product $L^\be\!\otimes_{\rm h}\!L^\be\!\otimes_{\rm h}\!L^\be$ for an arbitrary function $f$ in $B_{\be,1}^1(\R^2)$. 

\begin{thm}
\label{kpp<2}
Let $1\le p<2$.
There are spectral measures $E_1$, $E_2$ and $E_3$ on Borel subsets of $\R$, a function $\Phi$ in the Haagerup tensor product 
$L^\be(E_1)\!\otimes_{\rm h}\!L^\be(E_2)\!\otimes_{\rm h}\!L^\be(E_3)$ and an operator $Q$ in $\bS_p$ such that
$$
\iiint\Phi(x_1,x_2,x_2)\,dE_1(x_1)\,dE_2(x_2)Q\,dE_3(x_3)\not\in\bS_p.
$$
\end{thm}

\Pf Assume the contrary. Then the linear operator
$$
Q\mapsto\iiint\Phi(x_1,x_2,x_2)\,dE_1(x_1)\,dE_2(x_2)Q\,dE_3(x_3)
$$
is bounded on $\bS_p$ for arbitrary Borel spectral measures $E_1$, $E_2$, and $E_3$
and for an arbitrary function $\Phi$ in 
$L^\be(E_1)\!\otimes_{\rm h}\!L^\be(E_2)\!\otimes_{\rm h}\!L^\be(E_3)$. 
Suppose now that $\Psi$ belongs to the Haagerup-like tensor product
$L^\be(E_1)\!\otimes_{\rm h}\!L^\be(E_2)\!\otimes^{\rm h}\!L^\be(E_3)$
of the first kind. For a finite rank operator $T$ consider the triple operator integral
$$
W=\iint\!\!\upint\Psi(x_1,x_2,x_3)\,dE_1(x_1)T\,dE_2(x_2)\,dE_3(x_3).
$$
We define the function $\Phi$ defined by
$$
\Phi(x_2,x_3,x_1)=\Psi(x_1,x_2,x_3).
$$
Let $Q\in\bS_p$. We have
\begin{align*}
\trace(WQ)&=
\trace\left(\left(
\iiint
\Psi(x_1,x_2,x_3)\,dE_2(x_2)\,dE_3(x_3)Q\,dE_1(x_1)
\right)T\right)\\[.2cm]
&=\trace\left(\left(
\iiint
\Phi(x_2,x_3,x_1)\,dE_2(x_2)\,dE_3(x_3)Q\,dE_1(x_1)
\right)T\right)
\end{align*}
(see the definition of triple operator integrals with integrands in the 
Haagerup-like tensor product of the first kind in \S~\ref{Ttoi}).

Thus
\begin{align*}
|\trace(WQ)|&=
\left|\trace\left(\left(
\iiint
\Phi(x_2,x_3,x_1)\,dE_2(x_2)\,dE_3(x_3)Q\,dE_1(x_1)
\right)T\right)\right|\\[.2cm]
&\le
\left\|\left(
\iiint
\Phi(x_2,x_3,x_1)\,dE_2(x_2)\,dE_3(x_3)Q\,dE_1(x_1)
\right)\right\|_{\bS_p}\|T\|_{\bS_{p'}}\\[.2cm]
&\le\|\Phi\|_{L^\be\!\otimes_{\rm h}\!L^\be\!\otimes_{\rm h}\!L^\be}
\|Q\|_{\bS_p}\|T\|_{\bS_{p'}}
\end{align*}
(throughout the proof of this theorem in the case $p=1$, the norm in $\bS_{p'}$ has to be replaced with the operator norm).

It follows that
\begin{align}
\label{WSp'}
\|W\|_{\bS_{p'}}&=
\left\|
\iint\!\!\upint\Psi(x_1,x_2,x_3)\,dE_1(x_1)T\,dE_2(x_2)\,dE_3(x_3)
\right\|_{\bS_{p'}}\nonumber\\[.2cm]
&\le\|\Psi\|_{L^\be\!\otimes_{\rm h}\!L^\be\!\otimes^{\rm h}\!L^\be}\|T\|_{\bS_{p'}}.
\end{align}

By Theorem \ref{Bes}, 
$\dg^{[1]}f\in L^\be\!\otimes_{\rm h}\!L^\be\!\otimes^{\rm h}\!L^\be$ for every
$f$ in $B_{\be,1}^1(\R^2)$ and by \rf{perfor},
\begin{align*}
f(A_1,B)&-f(A_2,B)\nonumber\\[.2cm]
&=\iint\!\!\upint\big(\dg^{[1]}f\big)(x_1,x_2,y)
\,dE_{A_1}(x_1)(A_1-A_2)\,dE_{A_2}(x_2)\,dE_B(y)
\end{align*}
for arbitrary finite rank self-adjoint operators $A_1$, $A_2$, and $B$.
It remains to observe that by inequality \rf{WSp'},
\begin{align*}
\|f(A_1,B_1)-f(A_2,B)\|_{\bS_{p'}}
&\le
\|\dg^{[1]}f\|_{L^\be\!\otimes_{\rm h}\!L^\be\!\otimes^{\rm h}\!L^\be}
\|A_1-A_2\|_{\bS_{p'}}\\[.2cm]
&\le\const\|f\|_{B_{\be,1}^1}\|A_1-A_2\|_{\bS_{p'}}
\end{align*}
which contradicts Theorem \ref{BnLte}. $\bl$

If we pass to the adjoint operator, we can see that for $p\in[1,2)$, there exist a function $\Psi$ in the Haagerup tensor product 
$L^\be\!\otimes_{\rm h}\!L^\be\!\otimes_{\rm h}\!L^\be$ and an operator $Q$ in $\bS_p$ such that
$$
\iiint\Phi(x_1,x_2,x_2)\,dE_1(x_1)Q\,dE_2(x_2)\,dE_3(x_3)\not\in\bS_p.
$$

The following application of Theorem \ref{BnLte} shows that for functions
$f$ in $B_{\be,1}^1(\R^2)$, the divided differences $\dg^{[1]}f$ and $\dg^{[2]}f$
do not have to belong to the Haagerup tensor product 
$L^\be\!\otimes_{\rm h}\!L^\be\!\otimes_{\rm h}\!L^\be$. We state the result for
$\dg^{[1]}f$. 

\begin{thm}
\label{Bnd12}
There exists a function $f$ in the Besov class $B_{\be,1}^1(\R^2)$ such that the divided difference $\dg^{[1]}f$ does not belong to 
$L^\be\!\otimes_{\rm h}\!L^\be\!\otimes_{\rm h}\!L^\be$.
\end{thm}

\Pf Assume the contrary. Then the map
$$
f\mapsto\dg^{[1]}f
$$
is a bounded linear operator from $B_{\be,1}^1(\R^2)$ to 
$L^\be\!\otimes_{\rm h}\!L^\be\!\otimes_{\rm h}\!L^\be$.

By \rf{perfor},
\begin{align*}
f(A_1,B)&-f(A_2,B)\nonumber\\[.2cm]
&=\iint\!\!\upint\big(\dg^{[1]}f\big)(x_1,x_2,y)
\,dE_{A_1}(x_1)(A_1-A_2)\,dE_{A_2}(x_2)\,dE_B(y)
\end{align*}
for arbitrary finite rank self-adjoint operators $A_1$, $A_2$, and $B$.
It follows now from inequality \rf{opno} that
$$
\|f(A_1,B)-f(A_2,B)\|
\le\big\|\dg^{[1]}f\big\|_{L^\be\!\otimes_{\rm h}\!L^\be\!\otimes_{\rm h}\!L^\be}
\|A_1-A_2\|
\le\const\|f\|_{B_{\be,1}^1}\|A_1-A_2\|
$$
which contradicts Theorem \ref{BnLte}. $\bl$

\medskip

{\bf Remark.} It is easy to observe that the construction given in the proof of Theorem \ref{nLte} allows us to construct a function $f$ in $B_{\be,1}^1(\R^2)$, for which both divided differences $\dg^{[1]}f$ and $\dg^{[2]}f$
do not belong to the Haagerup tensor product 
$L^\be\!\otimes_{\rm h}\!L^\be\!\otimes_{\rm h}\!L^\be$.

\

\section{\bf Points of Lipschitzness}
\setcounter{equation}{0}
\label{tlL}

\

We have shown in \S~\ref{Bp>2} that for functions $f$ in $B_{\be,1}^1(\R^2)$ there is no Lipschitz type estimate in the operator norm. It turns out however that for certain pairs $(A_*,B_*)$ of self-adjoint operators the function
$(A,B)\mapsto f(A,B)$ is Lipschitz at $(A_*,B_*)$ for all functions 
$f$ in $B_{\be,1}^1(\R^2)$. We establish in this section the fact that the pairs 
$(\a I,\b I)$ are points of Lipschitzness for all $\a$ and $\b$ in $\R$. The same is true in the Schatten--von Neumann norm of $\bS_p$ (quasi-norm for $p<1$) for all $p>0$.

\begin{thm}
\label{ttlL}
There exists a positive number $C$ such that
$$
\|f(A,B)-f(\a I,\b I)\|\le C\|f\|_{B^1_{\be,1}}\max\big\{\|A-\a I\|,\|B-\b I\|\big\}
$$
for arbitrary $f$ in $B_{\be,1}^1(\R^2)$, for arbitrary self-adjoint operators $A$ and $B$, and for arbitrary real numbers $\a$ and $\b$.
\end{thm}

\begin{thm}
\label{tlLSp}
Let $f\in B^1_{\be1}(\R^2)$. Then there exists a positive number $C$ such that
$$
\|f(A,B)-f(\a I,\b I)\|_{\bS_p}\le C\|f\|_{B^1_{\be1}}
\max\big\{\|A-\a I\|_{\bS_p},\|B-\b I\|_{\bS_p}\big\}
$$
for arbitrary $f$ in $B_{\be,1}^1(\R^2)$, for arbitrary self-adjoint operators $A$ and $B$, and for arbitrary $p>0$, and $\a,\,\b\in\R$.
\end{thm}

First we obtain several auxiliary results.

\begin{lem} 
Let $f\in B^1_{\be1}(\R^2)$. Then 
$\big(\dg^{[1]}f\big)(0,\cdot,\cdot)\in L^\be(\R)\!\otimes_{\rm h}\!L^\be(\R)$ and
$$
\big\|\big(\dg^{[1]}f\big)(0,\cdot,\cdot)\big\|_{L^\be\!\otimes_{\rm h}\!L^\be}\le\const\|f\|_{B^1_{\be1}}.
$$
\end{lem}

\Pf It suffices to prove that
\bay
\label{75}
\big\|\big(\dg^{[1]}f\big)(0,\cdot,\cdot)\big\|_{L^\be\!\otimes_{\rm h}\!L^\be}\le\const\|f\|_{L^\be}
\ey
for an arbitrary bounded function $f$ with  $\supp\F f\subset\{\xi\in\R^2:\|\xi\|\le1\}$.
By Theorem \ref{rrp}, we have
$$
\big(\dg^{[1]}f\big)(0,x,y)=\big(\dg^{[1]}f\big)(x,0,y)=
\sum_{j\in\Z}\frac{\sin(x-j\pi)}{x-j\pi}\big(\dg^{[1]}f\big)(0,j\pi,y).
$$
Now \rf{75} follows from \rf{sinusy} because
$$
\sum_{j\in\Z}\big|\big(\dg^{[1]}f\big)(0,j\pi,y)\big|^2
\le\|f\|_{L^\be}^2+2\left(\frac4{\pi^2}\sum_{j=1}^\be\frac1{j^2}\right)\|f\|_{L^\be}^2
=\frac73\|f\|_{L^\be}^2.\quad\bl
$$

\begin{cor}  
Let $f\in B^1_{\be1}(\R^2)$ and let $\a,\,\b\in\R$. Then both
$\big(\dg^{[1]}f\big)(\a,\cdot,\cdot)$
and $\big(\dg^{[2]}f\big)(\cdot,\cdot,\b)$ belong to 
$L^\be(\R)\!\otimes_{\rm h}\!L^\be(\R)$. Moreover, there exists a positive number $C$ such that
$$
\big\|\big(\dg^{[1]}f\big)(a,\cdot,\cdot)\big\|_{L^\be\!\otimes_{\rm h}\!L^\be}\le C\|f\|_{B^1_{\be1}}\quad\text{and}
\quad\big\|\big(\dg^{[2]}f\big)(\cdot,\cdot,b)\big\|_{L^\be\otimes_{\rm h}L^\be}\le C\|f\|_{B^1_{\be1}}
$$
for all $a,b\in\R$.
\end{cor}

\begin{cor} 
For arbitrary self-adjoint operators $A$ and $B$,
$$
\|f(A,B)-f(\a I,B)\|\le C\|A-\a I\|\cdot\|f\|_{B^1_{\be1}}
$$
and
$$
\|f(A,B)-f(A,\b I)\|\le C\|B-\b I\|\cdot\|f\|_{B^1_{\be1}}.
$$
\end{cor}

\Pf Clearly, it suffices to prove the first inequality for $a=0$. Since
$$
f(x,y)-f(0,y)=x\,\dg^{[1]}f(0,x,y),\quad x,~y\in\R,
$$
it follows that
\bay
\label{fABdg}
f(A,B)-f(\0,B)=A(\dg^{[1]}f\big)(0,A,B).
\ey
Since functions in $L^\be(\R)\!\otimes_{\rm h}\!L^\be(\R)$ are Schur multipliers (see Subsection 2.3) we have
\begin{align}
\label{tsepo}
\|f(A,B)-f(\0,B)\|&\le\|A\|\cdot\|(\dg^{[1]}f\big)(0,A,B)\|\nonumber\\[.2cm]
&\le\|A\|\cdot\big\|(\dg^{[1]}f\big)(0,\cdot,\cdot)\big\|_{\fM(E_A,E_B)}\nonumber\\[.2cm]
&\le\|A\|\cdot\big\|(\dg^{[1]}f\big)(0,\cdot,\cdot)\big\|_{L^\be\otimes_{\rm h}L^\be}\le 
C\|A\|\cdot\|f\|_{B^1_{\be1}}.
\end{align}
The second inequality can be proved in the same way. $\bl$

It turns out that similar estimates hold in the norm of $\bS_p$ for arbitrary $p>0$.

\begin{cor} 
\label{2ner}
For every self-adjoint operators $A$ and $B$
$$
\|f(A,B)-f(\a I,B)\|_{\bS_p}\le C\|A-\a I\|_{\bS_p}\|f\|_{B^1_{\be1}}
$$
and
$$
\|f(A,B)-f(A,\b I)\|_{\bS_p}\le C\|B-\b I\|_{\bS_p}\|f\|_{B^1_{\be1}}
$$
for all $p\in(0,+\be)$.
\end{cor}

\Pf Again, without loss of generality we may assume that $\a=0$. By \rf{fABdg}, 
$$
\|f(A,B)-f(\0,B)\|_{\bS_p}\le C\|A\|_{\bS_p}\|(\dg^{[1]}f\big)(0,A,B)\|
\le\|A\|_{\bS_p}\|f\|_{B^1_{\be1}}
$$
which is a consequence of \rf{tsepo}. The second inequality can be proved in exactly the same way. $\bl$

{\bf Proof of Theorem \ref{ttlL}.} Clearly,
$$
f(A,B)-f(\a I,\b I)=(f(A,B)-f(\a I,B))+(f(\a I,B)-f(\a I,\b I)).
$$
The result follows now from Corollary \ref{2ner}. $\bl$

Theorem \ref{tlLSp} can be proved in exactly the same way.

\

\section{\bf A sufficient condition for Lipschitz type estimates}
\setcounter{equation}{0}
\label{dostu}

\

We have seen in \S~\ref{Bp>2} that for functions $f$ in the Besov class 
$B_{\be,1}^1(\R^2)$, there is no Lipschitz type estimate in the operator norm as well as in the norm of $\bS_p$ for $p>2$.
In this section we obtain a simple sufficient condition for Lipschitz type estimates of the form \rf{Lipp>2} to hold in $S_p$ for every $p\ge1$ and in the operator norm.

We define a function class ${\mathscr C}$. Note that a similar class was defined in \cite{Pe0}.

\medskip

{\bf Definition.} {\it The class ${\mathscr C}$ of function on $\R^2$ is defined by
$$
{\mathscr C}\df B_{\be,1}^1(\R)\hat\otimes L^\be(\R)
\cap L^\be(\R)\hat\otimes B_{\be,1}^1(\R).
$$
In other words a function $f$ belongs to ${\mathscr C}$ if there are sequences of functions $\f_n$, $\psi_n$, $\f_n^\sharp$, and $\psi_n^\sharp$ on $\R$ such that
\bay
\label{rav}
f(x,y)=\sum_n\f_n(x)\psi_n(y)=\sum_n\f^\sharp_n(x)\psi^\sharp_n(y),\quad(x,y)\in\R^2,
\ey
\bay
\label{summy}
\sum_n\|\f_n\|_{B_{\be,1}^1(\R)}\|\psi_n\|_{L^\be(\R)}+
\sum_n\|\f_n^\sharp\|_{L^\be(\R)}\|\psi_n^\sharp\|_{B_{\be,1}^1(\R)}<\be.
\ey
The norm $\|f\|_{\mathscr C}$ of $f$ in the space ${\mathscr C}$ is, by definition
the infimum of {\em\rf{summy}} over all functions $\f_n$, $\psi_n$, $\f_n^\sharp$, and $\psi_n^\sharp$ satisfying {\em\rf{rav}}.
}

\begin{thm}
\label{LipSp}
There exists a positive number $C$ such that  
\bay
\label{neropno}
\|f(A_1,B_1)-f(A_2,B_2)\|\le C\|f\|_{\mathscr C}(\|A_1-B_1\|+\|A_2-B_2\|),
\ey
whenever $f\in{\mathscr C}$ and $A_1$, $A_2$, $B_1$, and $B_2$ are self-adjoined operators.

If $p\ge1$, $A_1-A_2\in\bS_p$, and $B_1-B_2\in\bS_p$, then
\bay
\label{nerSvNno}
\|f(A_1,B_1)-f(A_2,B_2)\|_{\bS_p}
\le C\|f\|_{\mathscr C}(\|A_1-B_1\|_{\bS_p}+\|A_2-B_2\|_{\bS_p}).
\ey
\end{thm}

\Pf Suppose that $f\in{\mathscr C}$ and \rf{summy} holds. Clearly,
$$
\|f(A_1,B_1)-f(A_2,B_2)\|\le\|f(A_1,B_1)-f(A_2,B_1)\|+\|f(A_2,B_1)-f(A_2,B_2)\|.
$$
Making use of inequality \rf{BesSch}, we obtain
\begin{align*}
\|f(A_1,B_1)-f(A_2,B_1)\|&\le\sum_n\|\f_n(A_1)-\f_n(A_2)\|\cdot\|\psi_n(B_1)\|
\\[.2cm]
&\le\const\sum_n\|\f_n\|_{B_{\be,1}^1(\R)}\|\psi_n\|_{L^\be(\R)}\|A_1-A_2\|.
\end{align*}
Similarly,
$$
\|f(A_2,B_1)-f(A_2,B_2)\|
\le\const\sum_n\|\f_n^\sharp\|_{L^\be(\R)}\|\psi_n^\sharp\|_{B_{\be,1}^1(\R)}.
$$
This implies \rf{neropno}. The proof of \rf{nerSvNno} is exactly the same. $\bl$

\

\section{\bf Functions of noncommuting unitary operators}
\setcounter{equation}{0}
\label{unit}

\

In this section we briefly explain that analogs of the above results hold for functions of noncommuting unitary operators hold. 

Suppose that $f$ is a function on $\T^2$ that belongs to the Besov space 
$B_{\be,1}^1(\T^2)$ (see Subsection 2.1.2). As we have observed in Subsection 2.4,
$f$ is a Schur multiplier with respect to arbitrary spectral Borel measures on $\T$.
This allows us to define functions $f(U,V)$ for (not necessarily commuting) unitary operators $U$ and $V$ on Hilbert space by the formula
$$
f(U,V)=\int\limits_\T\int\limits_\T f(\z,\t),dE_U(\z)\,dE_V(\t),
$$
where $E_U$ and $E_V$ are the spectral measures of $U$ and $V$.

As in the case of functions of self-adjoint operators, we are would like to use the formula:
\begin{align}
\label{osnfuni}
f(U_1,V_1)&-f(U_2,V_2)\nonumber\\[.2cm]
&=
\iint\!\!\upint\big(\dg^{[1]}f\big)(\z_1,\z_2,\t)
\,dE_{U_1}(\z_1)(U_1-U_2)\,dE_{U_2}(\z_2)\,dE_{V_1}(\t),\nonumber\\[.2cm]
&+\upint\!\!\!\iint\big(\dg^{[2]}f\big)(\z,\t_1,\t_2)
\,dE_{U_2}(\z)\,dE_{V_1}(\t_1)(V_1-V_2)\,dE_{V_2}(\t_2),
\end{align}
where the divided differences $\dg^{[1]}f$ and $\dg^{[2]}f$ are the functions on 
$\T^2$ defined by
$$
\big(\dg^{[1]}f\big)(\z_1,\z_2,\t)\df\frac{f(\z_1,\t)-f(\z_2,\t)}{\z_1-\z_2}
\quad\mbox{and}\quad
\big(\dg^{[2]}f\big)(\z,\t_1,\t_2)\df\frac{f(\z,\t_1)-f(\z,\t_2)}{\t_1-\t_2}.
$$

To establish formula \rf{osnfuni}, we should prove that $\dg^{[1]}f$ 
belongs to the Haagerup-like tensor product 
$L^\be\!\otimes_{\rm h}\!L^\be\!\otimes^{\rm h}\!L^\be$ of the first kind
and $\dg^{[2]}f$ belongs to the Haagerup-like tensor product 
$L^\be\!\otimes^{\rm h}\!L^\be\!\otimes_{\rm h}\!L^\be$ of the second kind
with respect to arbitrary Borel spectral measures.

To this end, we introduce the functions $\Xi_n$ on $\T$ defined by
$$
\Xi_n(z)\df\frac{z^{n+1}-z^{-n}}{(2n+1)(z-1)}=\frac1{2n+1}\sum_{k=-n}^nz^k,\quad z\in\T.
$$
For a positive integer $k$, we denote by $\Pi_k$ the group of $k$th roots of 1:
$$
\Pi_k\df\{\z\in\C:~\z^k=1\}.
$$

\begin{thm}
\label{razperu}
Let $n$ be a positive integer and let $f$ be a bounded function on $\T^2$ whose Fourier transform is supported 
in $\{(j,k)\in\Z^2: |j|\le n\}$. Then
$$
\big(\dg^{[1]}f\big)(\z_1,\z_2,\t)=
\sum_{\vk,\xi\in\Pi_{2n+1}}\Xi_n(\z_1\ov\vk)\,\Xi_n(\z_2\ov\xi)
\frac{f(\vk,\t)-f(\xi,\t)}{\z-\xi}.
$$
Moreover,
\bey
\sum_{\vk\in\Pi_{2n+1}}|\Xi_n(\z_1\ov\vk)|^2
=\sum_{\xi\in\Pi_{2n+1}}|\Xi_n(\z_2\ov\xi)|^2=1,
\quad \z_1,~\z_2\in\T,
\eey
and
\bey
\sup_{w\in\T}\left\|\left\{\frac{f(\vk,\t)-f(\xi,\t)}{\vk-\xi}
\right\}_{\vk,\xi\in\Pi_{2n+1}}\right\|_{\mB}\le\const(2n+1)\|f\|_{L^\be(\T)}.
\eey
\end{thm}

By $\|\cdot\|_{\mB}$ we mean the operator norm in the space of $(2n+1)\times(2n+1)$ matrices.

Theorem \ref{razperu} can be proved straightforwardly. We leave it as an exercise.

\begin{cor}
\label{sleu}
Let $f$ be a trigonometric polynomial of degree at most $n$ in each variable.
Then
\bay
\label{tenner1}
\big\|\dg^{[1]}f\big\|_{L^\be\!\otimes_{\rm h}\!L^\be\!\otimes^{\rm h}\!L^\be}
\le\const n\|f\|_{L^\be(\T^2)}.
\ey
\end{cor}

Similarly, it can be shown that under the hypotheses of Corollary \rf{sleu} 
\bay
\label{tenner2}
\big\|\dg^{[2]}f\big\|_{L^\be\!\otimes^{\rm h}\!L^\be\!\otimes_{\rm h}\!L^\be}
\le\const n\|f\|_{L^\be(\T^2)}.
\ey

Inequalities \rf{tenner1} and \rf{tenner2} imply the following result:

\begin{thm}
\label{teordlaB}
Let $f$ be a function in $B_{\be,1}^1(\T^2)$. Then $\dg^{[1]}f\in L^\be\!\otimes_{\rm h}\!L^\be\!\otimes^{\rm h}\!L^\be$,
$\dg^{[2]}f\in L^\be\!\otimes^{\rm h}\!L^\be\!\otimes_{\rm h}\!L^\be$,
$$
\big\|\dg^{[1]}f\big\|_{L^\be\!\otimes_{\rm h}\!L^\be\!\otimes^{\rm h}\!L^\be}
\le\const \|f\|_{B_{\be,1}^1}\quad\mbox{and}\quad
\big\|\dg^{[2]}f\big\|_{L^\be\!\otimes^{\rm h}\!L^\be\!\otimes_{\rm h}\!L^\be}
\le\const \|f\|_{B_{\be,1}^1}.
$$
\end{thm}

Theorem \ref{teordlaB} implies the following result:

\begin{thm}
\label{p2unit}
Let $1\le p\le2$ and let $f\in B_{\be,1}^1(\T^2)$. Suppose that $U_1,\,V_1,U_2,\,V_2$ are unitary operators such that $U_1-U_2\in\bS_p$ and $V_1-V_2\in\bS_p$. Then formula {\em\rf{osnfuni}} holds and
$$
\|f(U_1,V_1)-f(U_2,V_2)\|_{\bS_p}\le\const\|f\|_{B_{\be,1}^1}
\max\big\{\|U_1-U_2\|_{\bS_p},\|V_1-V_2\|_{\bS_p}\big\}.
$$
\end{thm}

As in the case of self-adjoint operators there is no Lipschitz type inequality  in the operator norm and in the norm of $\bS_p$ with $p>2$ for arbitrary functions 
$f$ in $B_{\be,1}^1(\T^2)$. To prove this, we can easily adjust the proof of Theorem \ref{BnLte} to the case of unitary operators.

\

\

\footnotesize
\noindent
\begin{tabular}{p{5cm}p{4.5cm}p{4.6cm}}
A.B. Aleksandrov & F.L. Nazarov &  V.V. Peller \\
St.Petersburg Branch & Department of Mathematics & Department of Mathematics  \\
Steklov Institute of Mathematics & Kent State University & Michigan State University\\
Fontanka 27 & Kent, Ohio 44242  & East Lansing, Michigan 48824 \\
 191023 St-Petersburg & USA & USA\\
 Russia
\end{tabular}

\end{document}